\documentclass[aos,MSNbibl,nameyear,seceqn,dvips]{arximspdf}
\usepackage{graphicx}
\doi{10.1214/13-AOS1100} 
\volume{41}
\issue{3}
\pubyear{2013}
\firstpage{1204}
\lastpage{1231}

\makeatletter
\newcommand{\rrvert}{\vert}
\newcommand{\llvert}{\vert}

\newtheorem{theorem}{Theorem}
\newtheorem{proposition}[theorem]{Proposition}
\newtheorem{lemma}[theorem]{Lemma}
\newproclaim{example}{Example}

\makeatother

\begin{document}
\begin{frontmatter}

\title{Asymptotic power of sphericity tests for high-dimensional data}
\runtitle{Asymptotic power of sphericity tests}

\begin{aug}
\author[A]{\fnms{Alexei} \snm{Onatski}\corref{}\ead[label=e1]{ao319@cam.ac.uk}},
\author[B]{\fnms{Marcelo J.} \snm{Moreira}\thanksref{m2}\ead[label=e2]{mjmoreira@fgv.br}}
\and
\author[C]{\fnms{Marc} \snm{Hallin}\thanksref{m3}\ead[label=e3]{mhallin@ulb.ac.be}}
\runauthor{A. Onatski, M. J. Moreira and M. Hallin}
\affiliation{University of Cambridge, FGV/EPGE, and Universit\'e libre
de Bruxelles and Princeton University}
\address[A]{A. Onatski\\
Faculty of Economics\\
University of Cambridge\\
Sidgwick Avenue\\
Cambridge, CB3 9DD\\
United Kingdom\\
\printead{e1}}
\address[B]{M. J. Moreira\\
Escola de P\'{o}s-gradua\c{c}\~{a}o em Economia\\
Funa\c{c}\~{a}o Getulio Vargas (FGV/EPGE)\\
Praia de Botafogo, 190-Sala 1100\\
Rio de Janeiro-RJ 22250-900\\
Brazil\\
\printead{e2}}
\address[C]{M. Hallin\\
ECARES\\
Universit\'e libre de Bruxelles CP 114/04\\
50, avenue F.D. Roosevelt\\
B-1050 Bruxelles\\
Belgium\\
\printead{e3}} 
\end{aug}

\thankstext{m2}{Supported by CNPq and NSF Grant SES-0819761.}

\thankstext{m3}{Supported by the Sonderforschungsbereich ``Statistical modelling of
nonlinear dynamic processes'' (SFB 823) of the Deutsche
Forschungsgemeinschaft, the IAP research network grant P7/06 of the
Belgian Federal Government (Belgian Science Policy), and a Discovery
grant of the Australian Research Council.}

\received{\smonth{8} \syear{2011}}
\revised{\smonth{6} \syear{2012}}

%
\begin{abstract}
This paper studies the asymptotic power of tests of sphericity against
perturbations in a single unknown direction as both the dimensionality
of the data and the number of observations go to infinity. We establish
the convergence, under the null hypothesis and contiguous alternatives,
of the log ratio of the joint densities of the sample covariance
eigenvalues to a Gaussian process indexed by the norm of the
perturbation. When the perturbation norm is larger than the
\textit{phase transition threshold} studied in Baik, Ben~Arous and
P{\'e}ch{\'e} [\textit{Ann. Probab.} \textbf{33}
(\citeyear{BaiBenPec05}) 1643--1697] the limiting process is
degenerate, and discrimination between the null and the alternative is
asymptotically certain. When the norm is below the threshold, the
limiting process is nondegenerate, and the joint eigenvalue
densities under the null and alternative hypotheses are mutually
contiguous. Using the asymptotic theory of statistical experiments, we
obtain asymptotic power envelopes and derive the asymptotic power for
various sphericity tests in the contiguity region. In particular, we
show that the asymptotic power of the Tracy--Widom-type tests is
trivial (i.e., equals the asymptotic size), whereas that of the
eigenvalue-based likelihood ratio test is strictly larger than the
size, and close to the power envelope.\looseness=-1
\end{abstract}

%

\end{frontmatter}

\section{Introduction}\label{sec1}

Recently, there has been much interest in testing sphericity in a
high-dimensional setting. Various tests have been proposed and analyzed
in \citet{LedWol02}, \citet{Sri05}, \citet{BirDet05},
\citet{Sch06}, \citet{Baietal09}, \citet{FisSunGal10},
\citet{CheZhaZho10} and \citet{BerRig}. In many studies, a
distinct interesting alternative to the null of sphericity is the
existence of a low-dimensional structure or signal in the data.
Detecting such a structure has been the focus of recent studies in
various applied fields including population and medical genetics
[\citet{PatPriRei06}], econometrics [Onatski (\citeyear{Ona09},
\citeyear{Ona10})], wireless communication [\citet{Biaetal}],
chemometrics [\citet{KriNad08}] and signal processing
[\citet{PerWol10}].

Most of the existing sphericity tests are based on the eigenvalues of the
sample covariance matrix, which constitute the maximal invariant statistic
with respect to orthogonal transformations of the data. The asymptotic
power of such tests
depends on the asymptotic behavior of the sample covariance eigenvalues
under the alternative hypothesis. When the alternative is a rank-$k$
perturbation of the null, the corresponding population covariance
matrix is
proportional to a sum of the identity matrix and a matrix of rank~$k$.
\citet{Joh01} calls such a situation ``spiked
covariance.''

The asymptotic behavior of the sample covariance eigenvalues in
``spiked covariance'' models of increasing
dimension is well studied. Consider the simplest case, when $k=1$. If the
largest population covariance eigenvalue is above the
``phase transition'' threshold studied in \citet{BaiBenPec05},
then the largest sample covariance eigenvalue remains separated from the
rest of the eigenvalues, which are asymptotically ``packed together as
in the support of the Marchenko--Pastur density''
[\citet{BaiSil06}]. Since the largest eigenvalue separates from the
``bulk,'' it is easy to detect a signal.

If the largest population covariance eigenvalue is at or below the
threshold, the empirical distribution of the sample covariance
eigenvalues still converges to the Marchenko--Pastur distribution, but
the largest sample covariance eigenvalue now converges to the upper
boundary of its support, both under the null of sphericity and the
``spiked'' alternative [\citet{SilBai95} and
\citet{BaiSil06}]. Hence, the signal detection becomes
problematic. At the threshold, the null and the alternative hypotheses
lead to different asymptotic distributions for the centered and
normalized largest sample covariance eigenvalue [\citet{BloVir}
and \citet{Mo12}], which implies some asymptotic
detection power. However, below the threshold, the difference
disappears with the joint distribution of any finite number of the
centered and normalized largest sample covariance eigenvalues
converging to the multivariate Tracy--Widom law under both the null and
the alternative [\citet{Joh01}, \citet{BaiBenPec05},
\citet{ElK07} and \citet{FerPec09}].

This similarity in the asymptotic behavior of covariance eigenvalues
under the null and the alternative prompts \citet{NadEde08} and
\citet{NadSil10} to call the transition threshold ``the
fundamental asymptotic limit of sample-eigenvalue-based detection.''
They claim that no reliable signal detection is possible below that
limit in the asymptotic sense. This asymptotic impossibility is also
pointed out and discussed in several other recent studies, including
\citet{PatPriRei06}, \citet{Hoy08}, \citet{Nad08},
\citet{KriNad09} and \citet{PerWol10}.

In this paper, we analyze the capacity of statistical tests to detect a
one-dimensional signal with the corresponding population covariance
eigenvalue below the ``impossibility
threshold,'' showing that the terminology
``impossibility threshold'' is overly pessimistic. We establish
that the eigenvalue region below the threshold actually is the region of
mutual contiguity [in the sense of \citet{LeC60}] of the joint distributions
of the sample covariance eigenvalues under the null and under the
alternative. We obtain the limit in distribution of the log likelihood
ratio process inside
this contiguity region and derive the asymptotic power envelope for
sample-eigenvalue-based detection tests.

The power envelope is larger than size for local alternatives and
monotonically tends to one as the signal's population eigenvalue
approaches the threshold from below. Hence, the detection of a signal with
high asymptotic probability is quite possible even in cases where the
largest population covariance eigenvalue is smaller than the threshold,
especially
when the distance from the threshold remains small.

In the contiguity region, the log likelihood ratio is asymptotically
equivalent to a simple statistic related to the Stieltjes transform
of the empirical distribution of the sample covariance eigenvalues. The
reason the asymptotic behavior of this statistic differs under the null
and under the alternative despite the apparent similarity of eigenvalue
behaviors just mentioned is that it is not based merely on a contrast
between the largest and the rest of the eigenvalues. The information about
the presence of the signal exploited by this statistic is hidden in the
small deviations of the empirical distribution of the eigenvalues from its
Marchenko--Pastur limit.

Let us examine our setting and our results in more detail. Suppose that
data consist of $n$ independent observations of $p$-dimensional
real-valued vectors $%
X_{t}$ distributed according to the Gaussian law with mean zero and
covariance matrix $\sigma^{2} ( I_{p}+hvv^{\prime} ) $,
where $%
I_{p}$ is the $p$-dimensional identity matrix, $\sigma$ and $h$ are
scalars and $v$ is a $p$-dimensional vector with Euclidean norm one. We are
interested in the asymptotic power of the tests of the null
hypothesis $H_{0}\dvtx h=0$ against the alternative $H_{1}\dvtx h>0$ based on the
eigenvalues of the sample covariance matrix of the data when both $n$
and $p$
go to infinity. The vector $v$ is an unspecified nuisance parameter indicating
the direction of the perturbation of sphericity. In contrast to \citet{BerRig},
who study signal detection in a similar setting where the vector $v$ is
sparse, we do not constrain
$v$ in any way except normalizing its Euclidean norm to one.

We consider the cases of known and unknown $\sigma^{2}$. For the sake of
brevity, in the rest of this Introduction, we discuss only the case of
unknown $\sigma^{2}$, which, in practice, is also more relevant. Let
$\lambda_{j}$ be the $j$th largest sample
covariance eigenvalue, let $\mu_{j}=\lambda_{j}/ ( \lambda
_{1}+\cdots+\lambda_{p} ) $ be its normalized version and let $\mu
= ( \mu_{1},\ldots,\mu_{m-1} ) $, where $m=\min( n,p
) $.
We begin our analysis with a study of the asymptotic properties of the
likelihood ratio process $L ( h;\mu) $ defined as the ratio of
the density of $\mu$ when $h\neq0$ to that when $h=0$. We represent $%
L ( h;\mu) $ in the form of an integral over a contour in the
complex plane and use the Laplace approximation method and recent results
from the large random matrix theory to derive an asymptotic expansion
of $%
L ( h;\mu) $ as $p,n\rightarrow\infty$ so that
$p/n\rightarrow
c\in( 0,\infty) $, which we throughout abbreviate into
$p,n\rightarrow_{c} \infty$.

We show that, for any $\bar{h}$ such that $0<\bar{h}<\sqrt{c}$, $\ln
L (
h;\mu) $ converges in distribution under the null to a Gaussian
process $\mathcal{L}(h;\mu)$ on $h\in[ 0,\bar{h} ] $ with
\[
\mathrm{E} \bigl[ \mathcal{L}(h;\mu) \bigr]
=\tfrac{1}{4} \bigl[ \ln\bigl( 1-c^{-1}h^{2} \bigr) +c^{-1}h^{2} \bigr]
\]
and
\[
\operatorname{Cov} \bigl(
\mathcal{L}(h_{1};\mu),\mathcal{L}(h_{2};\mu)\bigr) =-
\tfrac{1}{2} \bigl[ \ln\bigl( 1-c^{-1}h_{1}h_{2} \bigr) +c^{-1}h_{1}h_{2}\bigr].
\]
By Le Cam's first lemma [see \citet{van98}, page 88], this implies
that the joint distributions of the normalized sample covariance
eigenvalues under the null and under the alternative are mutually
contiguous for any $h\in[0,\bar{h}]$. We also show that these joint
distributions are not mutually contiguous for any $h>\sqrt{c}$.

Since $\mathcal{L}(h;\mu)$, as a likelihood ratio process,
is not
of the LAN Gaussian shift type, local asymptotic normality does not hold, and
the asymptotic optimality analysis of tests of $H_{0}\dvtx h=0$
against $H_{1}\dvtx h>0$ is difficult. However, an asymptotic power envelope is
easy to construct using the Neyman--Pearson lemma along with Le Cam's
third lemma.
We show that, for tests of asymptotic size $\alpha$, the maximum
achievable power
against a specific alternative $h=h_{1}$ is $1-\Phi[ \Phi
^{-1} ( 1-\alpha) -\sqrt{-\frac{1}{2} ( \ln(
1-c^{-1}h_{1}^{2} ) +c^{-1}h_{1}^{2} ) } ] $, where~$\Phi
$, as usual,
denotes the standard normal distribution function.

Using our result on the limiting distribution of $\ln L ( h;\mu ) $ and
Le Cam's third lemma, we compute the asymptotic powers of several
previously proposed tests of sphericity and of the likelihood ratio
(LR) test based on $\mu$. We find that the power of the LR test comes
close to the asymptotic power envelope. The LR test outperforms the
test proposed by \citet{Joh71} and studied in
\citet{LedWol02}, as well as \citet{Sri05} and the test
proposed by \citet{Baietal09}. The asymptotic powers of the tests
based on the largest sample covariance eigenvalue, such as the tests
proposed by \citet{Bej05}, \citet{PatPriRei06},
\citet{KriNad09}, \citet{Ona09}, \citet{Biaetal} and
\citet{NadSil10}, equals the tests' asymptotic size for
alternatives in the contiguity region.

The rest of the paper is organized as follows. Section~\ref{sec2} provides a
representation
of the likelihood ratio in terms of a contour integral. Section~\ref{sec3} applies
Laplace's method to obtain an asymptotic approximation to the contour
integral. Section~\ref{sec4} uses that approximation to establish
the convergence of the log likelihood ratio process to a Gaussian process.
Section~\ref{sec5} provides an analysis of the asymptotic power of various
sphericity tests and derives the asymptotic power envelope. Section
\ref{sec6}
concludes. Proofs are given in the \hyperref[app]{Appendix}; the more technical ones
are relegated to the Supplementary Appendix [\citet{OnaMorHal}].

\section{Likelihood ratios as contour integrals}\label{sec2}

Let $X$ be a $p\times n$ matrix with i.i.d. real Gaussian $N (
0,\sigma
^{2} ( I_{p}+hvv^{\prime} ) ) $ columns. Let $\lambda
_{1}\geq\lambda_{2}\geq\cdots\geq\lambda_{p}$ be the ordered eigenvalues
of $\frac{1}{n}XX^{\prime}$ and let $\lambda= ( \lambda
_{1},\ldots,\lambda_{m} ) $, where $m=\min\{ n,p \} $.
Finally, let $\mu= ( \mu_{1},\ldots,\mu_{m-1} ) $, where $\mu
_{j}=\lambda_{j}/ ( \lambda_{1}+\cdots+\lambda_{p} ) $.

As explained in the \hyperref[sec1]{Introduction}, our goal is to study the asymptotic power
of the eigenvalue-based tests of $H_{0}\dvtx h=0$ against $H_{1}\dvtx h>0$. If
$\sigma
^{2}$ is known, the model is invariant with respect to orthogonal
transformations, and the
maximal invariant statistic is $\lambda$. Therefore, we consider tests
based on $\lambda$. If $\sigma^{2}$ is unknown (which, strictly speaking,
is what is meant by ``sphericity''), the
model is invariant with respect to orthogonal transformations and
multiplications by nonzero
scalars, and the maximal invariant is $\mu$. Hence, we consider tests based
on $\mu$. Note that the distribution of $\mu$ does not depend on
$\sigma
^{2}$, whereas if $\sigma^{2}$ is known, we can always normalize
$\lambda$
dividing it by $\sigma^{2}$. Therefore, in what follows, we will assume
that $\sigma^{2}=1$ without loss of generality.

Let us denote the joint density of $\lambda_{1},\ldots,\lambda_{m}$ as $%
p ( \lambda;h ) $ and that of $\mu_{1},\ldots,\break\mu_{m-1}$ as $%
p ( \mu;h ) $. The following proposition gives explicit formulas
for $p ( \lambda;h ) $ and $p ( \mu;h ) $.

\begin{proposition}\label{pr1}
Let $\mathcal{S} ( r )$ be the
$(r-1)$-dimensional unit sphere, and let $ ( \mathrm{d}x_{r} )$
be the invariant measure on $\mathcal{S} ( r ) $
normalized so that the total measure is one. Further, let $\Lambda=%
\operatorname{diag} ( \lambda_{1},\ldots,\lambda_{p} ) $ and
$M=\operatorname{diag} ( \mu_{1},\ldots,\mu_{p} ) $.
Then%
%
\begin{equation}
\label{alternativereal}
p ( \lambda;h ) = \frac{\gamma( n,p,\lambda) }{%
( 1+h ) ^{n/2}}\int_{\mathcal{S}(p)}e^{({n}/{2})({h}/({%
1+h}))x_{p}^{\prime}\Lambda x_{p}}
( \mathrm{d}x_{p} )
\end{equation}
and
\begin{eqnarray}
\label{alternativerealmu}\quad
p ( \mu;h ) &=& \frac{\delta( n,p,\mu) }{ (
1+h ) ^{n/2}}\nonumber\\[-8pt]\\[-8pt]
&&{}\times\int_{0}^{\infty}y^{({np-2})/{2}}e^{-{n}y/{2}}
\int_{\mathcal{S}(p)}e^{({n}/{2})({yh}/({1+h}))x_{p}^{\prime
}Mx_{p}} ( \mathrm{d}x_{p} )
\,\mathrm{d}y,\nonumber
\end{eqnarray}
where $\gamma( n,p,\lambda) $ and $\delta
( n,p,\mu) $ depend only on $n$ and $p$,
and on $\lambda$ and $\mu$, respectively.
\end{proposition}

The spherical integrals in (\ref{alternativereal}) and
(\ref{alternativerealmu}) can be represented in the form of a confluent
hypergeometric function ${}_{1}F_{1}$ of matrix argument
[\citet{Hil01}, page~4]. For example, for the integral in
(\ref{alternativereal}),
\[
\int_{\mathcal{S}(p)}e^{({n}/{2})({h}/({1+h}))%
x_{p}^{\prime}\Lambda x_{p}} ( \mathrm{d}x_{p} ) =
{}_{1}F_{1} \biggl( \frac{1}{2},\frac{p}{2};
\frac{n}{2}\frac
{h}{1+h}\Lambda\biggr).
\]
\citet{ButWoo02} develop Laplace approximations to functions
${}_{1}F_{1}$ but do not analyze the asymptotic behavior of the approximation
errors. The next lemma derives an alternative
representation of the spherical integrals in Proposition~\ref{pr1}. This
representation has the form of a contour integral of a single complex
variable, and our asymptotic analysis will be based on the Laplace
approximation to such an integral.

\begin{lemma}\label{le2}
Let $D=\operatorname{diag} ( d_{1},\ldots,d_{r} )
$, where $d_{j}$ are arbitrary complex numbers. Further, let $%
\mathcal{K}$ be a contour in the complex plane starting at
%
\begin{figure}

\includegraphics{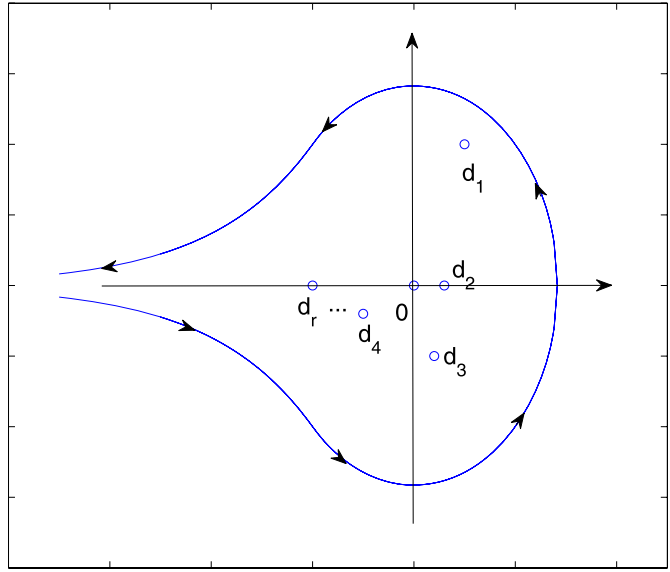}

\caption{Contour of integration $\mathcal{K}$ in
(\protect\ref{Erdelyirepresentation}).}\label{contourtheory}
\end{figure}
$-\infty$, encircling counter-clockwise the points
$0,d_{1},\ldots,d_{r}$,
and going back to $-\infty$. Such a contour is
shown in
Figure~\ref{contourtheory}. We have
%
\begin{equation}\label{Erdelyirepresentation}
\int_{\mathcal{S} ( r ) }e^{x_{r}^{\prime}Dx_{r}} ( \mathrm{d}%
x_{r} ) =\frac{\Gamma( r/2 ) }{2\pi i}\oint_{\mathcal{%
K}}e^{s}\prod
_{j=1}^{r} ( s-d_{j} )
^{-{1}/{2}}\,\mathrm{d}s.
\end{equation}
\end{lemma}

\begin{pf}
The integral on the left-hand side of (\ref{Erdelyirepresentation}) is the
expected value of $\exp( \frac{y_{1}^{2}d_{1}+\cdots+y_{r}^{2}d_{r}}{
y_{1}^{2}+\cdots+y_{r}^{2}} ) $, where $y_{1},\ldots,y_{r}$ are independent
standard normal random variables. The variables $u_{j}= \frac{y_{j}^{2}%
}{y_{1}^{2}+\cdots+y_{r}^{2}}$, $j=1,\ldots,r$, have Dirichlet distribution $%
\mathcal{D} ( k_{1},\ldots,k_{r}) $ with parameters $k_{1}=\cdots=k_{r}=\frac
{1}{2}%
$. Denoting the expectation operator with respect to such a
distribution as $%
\mathrm{E}_{\mathcal{D}}$, we have%
%
\begin{equation}\label{Dirichletaverage}
\int_{\mathcal{S} ( r ) }e^{x_{r}^{\prime}Dx_{r}} ( \mathrm{d}
x_{r} ) =\mathrm{E}_{\mathcal{D}} \exp( u_{1}d_{1}+\cdots+u_{r}d_{r}).
\end{equation}
Now, expanding the exponent in the latter expression into power series and
taking expectations term by term yields%
%
\begin{equation}\label{expexpansion}
\mathrm{E}_{\mathcal{D}} \exp( u_{1}d_{1}+\cdots+u_{r}d_{r}
) =\sum_{k=0}^{\infty}\frac{\mathrm{E}_{\mathcal{D}}
(u_{1}d_{1}+\cdots+u_{r}d_{r} ) ^{k}}{k!}.
\end{equation}
The Dirichlet average of $ ( u_{1}d_{1}+\cdots+u_{r}d_{r} ) ^{k}$ is
well studied. By Theorem~3.1 of \citet{Dic83},
%
\begin{eqnarray}\label{DT}
&&
\mathrm{E}_{\mathcal{D}} \bigl[ ( u_{1}d_{1}+\cdots+u_{r}d_{r})
^{k} \bigr]\nonumber\\[-8pt]\\[-8pt]
&&\qquad
= \mathop{\sum_{m_{1},\ldots,m_{r}\geq0}}_{m_{1}+\cdots+m_{r}=k}
\frac{k!}{%
m_{1}!\cdots m_{r}!}\frac{ ( 1/2 ) _{m_{1}}\cdots  ( 1/2 )
_{m_{r}}%
}{ ( r/2 ) _{k}}d_{1}^{m_{1}}\cdots d_{r}^{m_{r}},\nonumber
\end{eqnarray}
where $(k)_{s}=k ( k+1 )\cdots(k+s-1)$ is Pochhammer's notation for
the shifted factorial.

Combining (\ref{DT}) with (\ref{expexpansion}) and (\ref{Dirichletaverage}), we get%
%
\begin{eqnarray}\label{Lauricella}
\int_{\mathcal{S} ( r ) }e^{x_{r}^{\prime}Dx_{r}} ( \mathrm{d}%
x_{r} ) &=&\sum_{m_{1},\ldots,m_{r}\geq0}\frac{ ( 1/2 )
_{m_{1}}\cdots ( 1/2 ) _{m_{r}}}{ ( r/2 )
_{m_{1}+\cdots+m_{r}}}%
\frac{d_{1}^{m_{1}}\cdots d_{r}^{m_{r}}}{m_{1}!\cdots m_{r}!}
\nonumber\\[-8pt]\\[-8pt]
&=&{}_{r}\Phi( 1/2,\ldots,1/2;r/2;d_{1},\ldots,d_{r}),\nonumber
\end{eqnarray}
where the last equality is the definition of the confluent form of the
Lauricella $F_{D}$ function, denoted as ${}_{r}\Phi(\cdot)$. The
functions ${}_{r}\Phi( \cdot)$ were introduced by \citet{Erd37}
and are discussed by \citet{SriKar85}. In probability and
statistics, they were recently used to study the mean of a Dirichlet
process [see \citet{LijReg04} and references therein].

\citet{Erd37}, formula (8,6), establishes the following contour
integral representation of ${}_{r}\Phi( \cdot)$:
%
\begin{eqnarray}\label{Erdelyi}
&&
{}_{r}\Phi( k_{1},\ldots,k_{r};t;d_{1},\ldots,d_{r}
)\nonumber\\[-8pt]\\[-8pt]
&&\qquad= \frac{\Gamma( t ) }{2\pi i}\oint_{\mathcal{K}%
}e^{s}s^{-t+k_{1}+\cdots+k_{r}}\prod
_{j=1}^{r} ( s-d_{j} )
^{-k_{j}}\,\mathrm{d}s.\nonumber
\end{eqnarray}
Lemma~\ref{le2} follows from equalities (\ref{Lauricella}) and
(\ref{Erdelyi}).\vadjust{\goodbreak}
\end{pf}

The contour integral representation given in Lemma~\ref{le2} has been derived
independently
by \citet{Mo12} and \citet{Wan12}, who use it to study
the largest sample covariance eigenvalue
when the corresponding population eigenvalue equals
the critical threshold or lies above it. Our proof effectively
takes advantage of old results of
\citet{Dic83} and \citet{Erd37}, and thus is different from the proofs in
the above mentioned papers.

Using Lemma~\ref{le2} and Proposition~\ref{pr1}, we derive contour integral representations
for the likelihood ratios $L ( h;\lambda) =
p ( \lambda;h ) /p ( \lambda;0 ) $ and $L (
h;\mu
) = p ( \mu;h ) / p ( \mu;0 ) $. The quantity $
L ( h;\lambda) $ is the likelihood ratio based on $\lambda$ as
opposed to the entire data $X$. Similarly, $L ( h;\mu) $ is the
likelihood ratio based on $\mu$.

\begin{lemma}\label{le3}
Let $\mathcal{K}$ be a contour in the complex plane that starts at
$-\infty$, then encircles counter-clockwise the sample covariance
eigenvalues $\lambda_{1},\ldots,\lambda_{p}$, and goes back to $-\infty$.
In addition, we require that for any $z\in\mathcal{K}$,
$\operatorname{Re}z<\frac{1+h}{h}S$, where $\operatorname{Re}z$ denotes
the real part of $z\in\mathbb{C}$ and $S=\lambda
_{1}+\cdots+\lambda_{p}$. Then,%
%
\begin{equation}\label{LRcontourr1}\quad
L ( h;\lambda) =k_{1} \biggl( \frac{2}{n} \biggr)
^{({p-2})/{2}}%
\frac{1}{2\pi i}\oint_{\mathcal{K}}e^{({n}/{2})({h}/({1+h}))%
z}
\prod_{j=1}^{p} ( z-\lambda_{j} )
^{-{1}/{2}}\,\mathrm{d}z
\end{equation}
and
\begin{eqnarray}
\label{LRcontourr2}\qquad
L ( h;\mu) &=&k_{2}\frac{S^{({p-2})/{2}}}{2\pi i}\nonumber\\[-8pt]\\[-8pt]
&&{}\times\oint_{%
\mathcal{K}}e^{-(({np-p+2})/{2})\ln( 1-({h}/({1+h}))
({z}/{S}) )
}
\prod_{j=1}^{p} ( z-\lambda_{j} )
^{-{1}/{2}}\,\mathrm{d}%
z,\nonumber
\end{eqnarray}
where $k_{1}=h^{-({p-2})/{2}} ( 1+h ) ^{
({p-n-2})/{2}%
}\Gamma( p/2 ) $ and $k_{2}=k_{1}\frac{\Gamma(
( np-p+2 ) /2 ) }{\Gamma( np/2 ) }$.
\end{lemma}

Close inspection of the proof of Lemma~\ref{le3} reveals
that the right-hand side of (\ref{LRcontourr2}) depends on $\lambda
$ only
through $\mu$. Although it is possible to express $L ( h;\mu
) $
as an explicit function of $\mu$, the implicit form given in (\ref
{LRcontourr2}) is convenient because it allows us to use similar
methods for
the asymptotic analysis of the two likelihood ratios.

In the next two sections, we perform an asymptotic analysis of $L (
h;\lambda) $ and $L ( h;\mu) $ that relies on the Laplace
approximation of the contour integrals in Lem\-ma~\ref{le3} after those contours
have been
suitably deformed without changing the value of the integrals.

\section{Laplace approximation}\label{sec3}

In this section, we derive the Laplace approximations to the contour
integrals in Lemma~\ref{le3}. Laplace's method for contour integrals is discussed,\vadjust{\goodbreak}
for example, in Chapter 4 of \citet{Olv97}. The method describes an
asymptotic approximation to a contour integral $\oint_{\mathcal{K}%
}e^{-nf(z)}g(z)\,\mathrm{d}z$ as $n\rightarrow\infty$, where $f(z)$ and
$%
g(z) $ are analytic functions of $z$. The approximation is usually
based on
the part of the contour integral coming from a neighborhood of some
point $%
z_{0}\in\mathcal{K}$, where $z_{0}$ is such that $\frac{\mathrm
{d}}{\mathrm{%
d}z}f(z_{0})=0$ and $\operatorname{Re}f(z_{0})=\min_{z\in
\mathcal{K}}\operatorname{Re}f(z)$. For such a point to exist, one
might need to deform the contour so that, by Cauchy's theorem, the
value of the integral does not change. Typically, the deformation is
chosen so that $\operatorname{Re} ( -f(z) ) $ declines in the fastest
way possible as $z$ goes away from $z_{0}$ along the contour. For this
reason, the method is called the \textit{method of steepest descent}.

The contour integrals in (\ref{LRcontourr1}) and (\ref{LRcontourr2})
can be represented in the Laplace form with a deterministic function $f(z)$
and a random function $g(z)$ that converges to a log-normal random process
on the contour as $p,n\rightarrow_{c} \infty$. To see this, note that the
logarithm of the multiple product in (\ref{LRcontourr1}) and (\ref
{LRcontourr2}) equals $-\frac{1}{2}\sum_{j=1}^{p}\ln( z-\lambda
_{j} ) $. For each $z$, this expression is a special form of the
linear spectral statistic $\sum_{j=1}^{p}\varphi(\lambda_{j})$
studied by
\citet{BaiSil04}. According to the central limit theorem (Theorem
1.1) established in that paper, the random variable
%
\begin{equation}\label{deltadefinition}
\Delta_{p} ( z ) = \sum_{j=1}^{p}
\ln( z-\lambda_{j} ) -p\int\ln( z-\lambda) \,\mathrm{d}\mathcal{F}_{p}
( \lambda)
\end{equation}
converges in distribution to a normal random variable when
$p,n\rightarrow_{c}
\infty$. Here $\mathcal{F}_{p} ( \lambda) $ is the cumulative
distribution function of the Marchenko--Pastur distribution with a mass
of $\max(
0,1-c_{p}^{-1} ) $ at zero and
density%
%
\begin{equation}\label{Marchenko-Pastur}
\psi_{p} ( x ) =\frac{1}{2\pi c_{p}x}\sqrt{ ( b_{p}-x ) (
x-a_{p} ) },
\end{equation}
where $c_{p}=p/n$, $a_{p}= ( 1-\sqrt{c_{p}} ) ^{2}$ and $b_{p}=
( 1+\sqrt{c_{p}} ) ^{2}$.

Such a convergence suggests the following choices of $f(z)$ and $g(z)$ in
the Laplace forms of the integrals in (\ref{LRcontourr1}) and (\ref
{LRcontourr2}):%
%
\begin{equation}\label{fdefinition}
f(z)=-\frac{1}{2} \biggl( \frac{h}{1+h}z-c_{p}\int\ln( z-
\lambda) \,\mathrm{d}\mathcal{F}_{p} ( \lambda) \biggr)
\end{equation}
and
%
\begin{equation}
\label{gdefinition}\qquad g(z) = \cases{\displaystyle  \exp\biggl\{ -\frac{1}{2}
\Delta_{p} ( z ) \biggr\}, \qquad \mbox{for (\ref{LRcontourr1})},
\vspace*{2pt}\cr
\displaystyle \exp
\biggl\{ -\frac{np-p+2}{2}\ln\biggl( 1-\frac{h}{1+h}\frac
{z}{S}
\biggr) -\frac{n}{2}\frac{h}{1+h}z-\frac{1}{2}
\Delta_{p} ( z ) \biggr\}, \vspace*{3pt}\cr
\qquad\hspace*{77.5pt}\mbox{for (\ref{LRcontourr2})}.}
\end{equation}

As mentioned above, a particularly useful deformation of $\mathcal{K}$
passes through the point $z=z_{0}(h)$ where $\frac{\mathrm{d}}{\mathrm
{d}z}%
f ( z ) =0$. Taking the derivative of the\vadjust{\goodbreak} right-hand side of
(\ref{fdefinition}), we see that $z_{0}(h)$ must satisfy%
%
\begin{equation}\label{criticalpoint}
\frac{h}{1+h}+c_{p}m_{p}\bigl(z_{0}(h)
\bigr)=0,
\end{equation}
where $m_{p} ( z ) = \int\frac{1}{\lambda-z}\,\mathrm{d}%
\mathcal{F}_{p} ( \lambda) $ is the Stieltjes transform of the
Marchenko--Pastur distribution with parameter $c_{p}$. The properties
of $%
m_{p}(z)$ are well studied. In particular, the analytic expression for $
m_{p} ( z ) $ is known; see, for example, equation (2.3) in \citet{Bai93}. For $z\neq0$, which lies outside the support of $\mathcal{F}%
_{p} ( \lambda) $, we have%
%
\begin{equation}\label{Stijeltjesanalytic}
m_{p} ( z ) =\frac{-z-c_{p}+1+\sqrt{ ( z-c_{p}-1 )
^{2}-4c_{p}}}{2c_{p}z},
\end{equation}
where the branch of the square root is chosen so that the real and the
imaginary parts of $\sqrt{ ( z-c_{p}-1 ) ^{2}-4c_{p}}$ have the
same signs as the real and the imaginary parts of $z-c_{p}-1$, respectively.

Substituting (\ref{Stijeltjesanalytic}) into (\ref{criticalpoint}) and
solving for $z_{0}(h)$ when $h\!\in\!( 0,\sqrt{c_{p}})$, we~get\looseness=-1%
%
\begin{equation}\label{zzero}
z_{0}(h)=\frac{ ( 1+h ) ( c_{p}+h ) }{h}.
\end{equation}\looseness=0
When $h\geq\sqrt{c_{p}}$, there are no solutions to (\ref{criticalpoint})
that lie outside the support of $\mathcal{F}_{p} ( \lambda) $.
When $h=\sqrt{c_{p}}$, the right-hand side of (\ref{zzero}) equals
$(1+\sqrt{c_{p}}) ^{2}$, which lies exactly on the boundary of
the support of $\mathcal{F}_{p} ( \lambda) $. When $h>\sqrt{c_{p}}%
$, (\ref{zzero}) provides a solution to (\ref{criticalpoint}) only when
the branch of the square root in (\ref{Stijeltjesanalytic}) is chosen
differently. As can be verified using (\ref{fdefinition}) and
(\ref{Stijeltjesanalytic}), in such a case,
$\frac{\mathrm{d}}{\mathrm{d}z}f(z)$ is strictly negative at $z=z_{0} (
h ) $ given by (\ref{zzero}).

\begin{figure}

\includegraphics{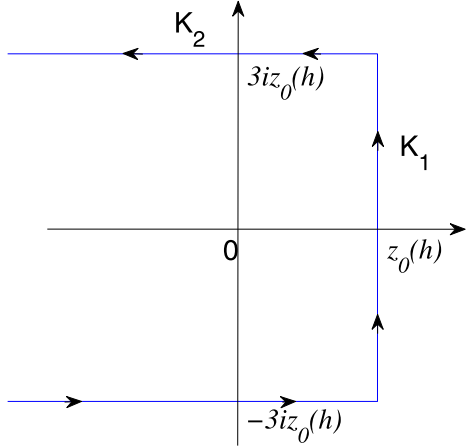}

\caption{Deformation $K$ of contour $\mathcal{K}$.}\label{contourdistortion}
\end{figure}

As $c_{p}\rightarrow c$, any fixed $h$ that is smaller than $\sqrt{c}$
eventually satisfies the inequality $h<\sqrt{c_{p}}$, so that $\frac
{\mathrm{d}}{%
\mathrm{d}z}f(z)=0$ at $z=z_{0}(h)$. Therefore, for $h<\sqrt{c}$, we will
deform the contour $\mathcal{K}$ into a contour $K$ that passes through
$%
z_{0}(h)$. We define $K$ as $K=K_{+}\cup K_{-}$, where $K_{-}$ is the
complex conjugate of $K_{+}$ and $K_{+}=K_{1}\cup K_{2}$ with%
%
\begin{equation}
\label{contours}
K_{1} = \bigl\{ z_{0}(h)+it\dvtx 0\leq t
\leq3z_{0}(h) \bigr\}
\end{equation}
and
\begin{equation}
\label{contours1}
K_{2} = \bigl\{ x+3iz_{0}(h)\dvtx -\infty<x
\leq z_{0}(h) \bigr\}.
\end{equation}
Figure~\ref{contourdistortion} illustrates the choice of $K$.

A proof of the following technical lemma is relegated to the Supplementary
Appendix [\citet{OnaMorHal}].\vspace*{-2pt}

\begin{lemma}\label{le4}
Suppose that our null hypothesis is true, and let $\bar
{h}$ be any fixed number such that $0<\bar{h}<\sqrt{c}$. Deforming
contour $\mathcal{K}$ into $K$ leaves the value of the
integrals (\ref{LRcontourr1}) and (\ref{LRcontourr2}) in Lemma
\ref{le3}
unchanged for all $h\in( 0,\bar{h} ] $
with probability approaching one as $p,n\rightarrow
_{c} \infty$.\vspace*{-2pt}
\end{lemma}

We now derive, uniform (over $h\in( 0,\bar{h} ] $), Laplace
approximations to the integrals (\ref{LRcontourr1}) and (\ref
{LRcontourr2}) in Lemma~\ref{le3}.\vadjust{\goodbreak} First, we introduce additional
notation. When $f ( z ) $ and $g(z)$ are analytic at $z_{0}=
z_{0} ( h ) $, let $f_{s}$ and $g_{s}$ with $s=0,1,\ldots$ be the
coefficients in the power series representations%
%
\begin{equation}\label{seriesforfandg}
f ( z ) =\sum_{s=0}^{\infty}f_{s} (
z-z_{0} ) ^{s},\qquad g(z)=\sum
_{s=0}^{\infty}g_{s} ( z-z_{0} )
^{s}.
\end{equation}
When $f ( z ) $ and $g(z)$ are not analytic at $z_{0}$, let the
coefficients $f_{s}$ and $g_{s}$ be arbitrary numbers for all $s$.

The following lemma is a generalization of the well-known Watson lemma for
contour integrals; see \citet{Olv97}, page 118. Theorem 7.1 in \citet{Olv97},
page~127, derives a similar generalization for the case when $f (
z ) $
and $g ( z ) $ are fixed deterministic analytic functions. In
contrast to
Olver's theorem, our lemma allows
$g(z)$ to be a random function, and $f(z)$ to depend on parameter $h$,
and obtains a uniform approximation over $h\in( 0,\bar{h} ]
$. The
proof is relegated to the Supplementary Appendix [\citet{OnaMorHal}].

\begin{lemma}\label{le5}
Under the conditions of Lemma~\ref{le4}, for any $h\in(
0,\bar{h}%
] $ and any positive integer $m$, as $p,n\rightarrow_{c} \infty$,
we have
%
\begin{equation}\label{Watson}\quad
\oint_{K}e^{-nf(z)}g(z)\,\mathrm{d}z=2e^{-nf_{0}} \Biggl[
\sum_{s=0}^{m-1}\Gamma\biggl( s+
\frac{1}{2} \biggr) \frac
{a_{2s}}{n^{s+1/2}}+%
\frac{O_{p} ( 1 ) }{hn^{m+1/2}}
\Biggr],
\end{equation}
where $O_{p} ( 1 ) $ is uniform in $h\in
(
0,\bar{h} ] $. The coefficients $a_{s}$ in
(\ref{Watson}) can be expressed through $f_{s}$ and $g_{s}$
defined above. In particular, we have
%
\begin{equation}\label{twocoefficients}
a_{0}=\frac{g_{0}}{2f_{2}^{1/2}}\quad\mbox{and}\quad a_{2}= \biggl
\{ 4g_{2}-%
\frac{6f_{3}g_{1}}{f_{2}}+ \biggl( \frac{15f_{3}^{2}}{2f_{2}^{2}}-
\frac
{6f_{4}%
}{f_{2}} \biggr) g_{0} \biggr\} \frac{1}{8f_{2}^{3/2}}.\hspace*{-28pt}\vadjust{\goodbreak}
\end{equation}
\end{lemma}

As we explained above, $z_{0}(h)$ is not a critical point of $f(z)$
when $h>%
\sqrt{c_{p}}$. This leads to a situation where the Laplace method for the
integral $\oint_{K}e^{-nf(z)}g(z)\,\mathrm{d}z$ delivers a
rather crude
approximation. Fortunately, our asymptotic analysis tolerates crude
approximations when $h>\sqrt{c_{p}}$. The following lemma, which is proven
in the Supplementary Appendix [\citet{OnaMorHal}], is sufficient for our purposes.

\begin{lemma}\label{le6}
Let $\tilde{h}>\sqrt{c}$, and denote by $K(\tilde {h})$ the
corresponding contour, as defined in (\ref{contours}) and
(\ref{contours1}). Under the null hypothesis, deforming the contour $
\mathcal{K}$ into $K(\tilde{h})$ leaves the value of the integrals in
Lemma~\ref{le3} unchanged for all $h\in[ \tilde{h},\infty
) $ with probability approaching one as $%
p,n\rightarrow_{c} \infty$. Further, for any $h\in[
\tilde{h},\infty
) $,%
%
\begin{equation}\label{Watson1}
\oint_{K(\tilde{h})}e^{-nf(z)}g(z)\,\mathrm{d}z=e^{-nf ( z_{0} (
\tilde{h} ) ) }O_{p}
( 1 ),
\end{equation}
where $O_{p} ( 1 ) $ is uniform over $h\in
[\tilde{h},\infty) $.
\end{lemma}

Neither Lemma~\ref{le5} nor Lemma~\ref{le6} addresses interesting cases with $h$ in a
neighborhood of $\sqrt{c}$. In such cases, $z_{0}(h)$ would be close to
the upper boundary of the support of the Marchenko--Pastur
distribution. This
may lead to the nonanalyticity of $f(z)$ and $g(z)$ on $K$ and a more
complicated asymptotic behavior of $g(z)$. We leave the analysis of cases
where $h$ may approach $\sqrt{c}$ for future research.

\citet{GuiMad05} study the asymptotic behavior of spherical
integrals using large deviation techniques. Their Theorems 3 and 6
imply Lemma~\ref{le6} and can be used to obtain the first term in the
asymptotic expansion of Lemma~\ref{le5}.

\section{Asymptotic behavior of the likelihood ratios}\label{sec4}

In this section, we discuss the asymptotic behavior of the likelihood
ratios $L ( h;\lambda) $ and $L ( h;\mu ) $. First, let us focus on the
case where $h\leq\bar{h}$. In the \hyperref[app]{Appendix}, we use
Lemmas~\ref{le4} and~\ref{le5} to derive the following theorem.

\begin{theorem}\label{th7}
Suppose that the null hypothesis is true ($h=0$). Let $\bar{h}$ be any
fixed number such that $0<\bar{h}<\sqrt{c}$ and let $C [
0,\bar{h} ] $ be the space of real-valued continuous
functions on $ [ 0,\bar{h} ] $
equipped with the supremum norm. Then as $%
p$, \mbox{$n\rightarrow_{c} \infty$},
we have, uniformly in $h\in(0, \bar{h}]$
%
\begin{equation}
\label{equivalence1}
L ( h;\lambda) =e^{-[ \Delta_{p} (
z_{0}(h) ) -\ln( 1-{h^{2}}/{c_{p}} ) ]/2
}+O_{p} \bigl( n^{-1}
\bigr)
\end{equation}
and
\begin{eqnarray}
\label{equivalence2}\quad
L ( h;\mu) &=&e^{-[ \Delta_{p} (
z_{0}(h) ) -\ln( 1-{h^{2}}/{c_{p}} )-
{h^{2}}/({2c_{p}%
})+({h}/{c_{p}}) ( S-p ) ]/2 }+O_{p} \bigl( n^{-1} \bigr).
\end{eqnarray}
Furthermore, $\ln L ( h;\lambda) $ and $\ln
L ( h;\mu) $, viewed as random elements of $C
[ 0,%
\bar{h} ], $ converge weakly to $\mathcal{L} (
h;\lambda) $ and $\mathcal{L} ( h;\mu) $
with Gaussian finite-dimensional\vadjust{\goodbreak} distributions such that, for any
$h_{1},\ldots,h_{r}\in[ 0,\bar{h} ] $,%
%
\begin{eqnarray}
\label{mean}
\mathrm{E} \bigl( \mathcal{L} ( h_{j};\lambda) \bigr) &=&
\tfrac{1%
}{4}\ln\bigl( 1-c^{-1}h_{j}^{2}
\bigr),
\\
\label{covariance}
\operatorname{Cov} \bigl( \mathcal{L} ( h_{j};\lambda),\mathcal
{L}%
( h_{k};\lambda) \bigr) &=&-\tfrac{1}{2}\ln\bigl(
1-c^{-1}h_{j}h_{k} \bigr),
\\
\label{meanmu}
\mathrm{E} \bigl( \mathcal{L} ( h_{j};\mu) \bigr) &=&
\tfrac
{1}{4}%
\bigl[ \ln\bigl( 1-c^{-1}h_{j}^{2}
\bigr) +c^{-1}h_{j}^{2} \bigr]
\end{eqnarray}
and
\begin{eqnarray}
\label{covariancemu}
\operatorname{Cov} \bigl( \mathcal{L} ( h_{j};\mu),\mathcal{L} (
h_{k};\mu) \bigr) &=&-\tfrac{1}{2} \bigl[ \ln\bigl(
1-c^{-1}h_{j}h_{k} \bigr) +c^{-1}h_{j}h_{k}
\bigr].
\end{eqnarray}
\end{theorem}

The log likelihood ratio processes studied in Theorem~\ref{th7} are not of the
\mbox{standard} \textit{locally asymptotically normal} form. This is because
they cannot be represented as $\varphi_{1}(h)W+\varphi_{2}(h)$, where $\varphi
_{1}(h)$ and $\varphi_{2}(h)$ are some deterministic functions of $h$,
and $%
W$ is a standard normal random variable. Indeed, had the representation
$%
\varphi_{1}(h)W+\varphi_{2}(h)$ been possible, the covariance of the
limiting log likelihood process at $h_{1}$ and $h_{2}$ would have been $
\varphi_{1}(h_{1})\varphi_{1}(h_{2})$. Hence, for $\mathcal{L} (
h;\lambda)$, for instance, we would have had $\varphi
_{1}(h)=\sqrt{-\frac{1}{%
2}\ln( 1-c^{-1}h^{2} ) }$ and $\varphi_{1}(h_{1})\varphi
_{1}(h_{2})=-\frac{1}{2}\ln( 1-c^{-1}h_{1}h_{2} ) $, which cannot
be true for all $0<h_{1}<\sqrt{c}$ and $0<h_{2}<\sqrt{c}$.

The quantity $\Delta_{p} ( z_{0}(h) ) $ plays an important
role in
the limits of experiments. The likelihood ratio processes are well
approximated by simple functions of $\Delta_{p} (
z_{0}(h) ) $ and $S$, which are easy to compute from the data and
are asymptotically Gaussian by the central limit theorem of \citet{BaiSil04}. Recalling the definition (\ref{deltadefinition})
of $\Delta_{p} (
z_{0}(h) ) $, we see that asymptotically, all statistical information
about parameter $h$ is contained in the deviations of the sample covariance
eigenvalues $\lambda_{1},\ldots,\lambda_{p}$ from $\lim_{n,p\rightarrow
\infty}z_{0}(h)=\frac{ ( 1+h ) ( h+c ) }{h}$. Although
the latter limit does not have an obvious interpretation when $h<\sqrt{c}$,
it is the probability limit of $\lambda_{1}$ under alternatives with
$h>%
\sqrt{c}$; see, for example, \citet{BaiSil06}.

Let us now consider cases where $h>\sqrt[.]{c}$. We prove the following
theorem in the \hyperref[app]{Appendix}.

\begin{theorem}\label{th8}
Suppose that the null hypothesis is true ($h=0$), and let $H$ be any
fixed number such that $\sqrt{c}<H<\infty$. Then as $p,n\rightarrow_{c}
\infty$, the following holds. For any $h\in[ H,\infty) $, the
likelihood ratios $L(h;\lambda) $ and $L ( h;\mu) $ converge to zero;
more precisely, there exists $\delta>0$ that depends
only on $H$ such that%
%
\begin{equation}\label{exponentialdecay}
L ( h;\lambda) =O_{p} \bigl( e^{-n\delta} \bigr)
\quad\mbox{and}\quad L ( h;\mu) =O_{p} \bigl( e^{-n\delta}
\bigr).
\end{equation}
\end{theorem}

Note that Theorem~\ref{th7} and Le Cam's first lemma [see \citet{van98},
page 88] imply that the joint distributions of $\lambda
_{1},\ldots,\lambda_{m}$
(as well as those of $\mu_{1},\ldots,\mu_{m-1}$)\vadjust{\goodbreak} under the null and
under the
alternative are mutually contiguous for any $h\in[ 0,\sqrt{c} ) $.
In contrast, Theorem~\ref{th8} shows that mutual contiguity is lost for $h>\sqrt{c}$.
For such $h$, consistent tests (as $p,n\rightarrow_{c} \infty$) exist
at any probability level $\alpha>0$.

In a similar setting, \citet{NadEde08} call the number of
``signal eigenvalues'' of the population
covariance matrix that exceed $1+\sqrt{c}$ the ``effective
number of identifiable signals'' [see also \citet{NadSil10}]. Theorems~\ref{th7} and~\ref{th8} shed light on the formal statistical
content of this concept. The ``identifiable
signals'' are detected with probability approaching
one in
large samples (irrespective of the probability level $\alpha>0$ at which
identification tests are performed). Other signals still can be detected,
but the probability of detecting them will never approach one (whatever the
probability level $\alpha<1$).

\section{Asymptotic power analysis}\label{sec5}

Theorem~\ref{th7} can be used to study ``local''
powers of the tests for detecting signals in noise. The nonstandard
form of
the limit of log likelihood ratio processes in our setting makes it hard
to develop tests with optimal local power properties. However, using the
Neyman--Pearson lemma and Le Cam's third lemma, we can analytically derive
the local asymptotic power envelope and compare local asymptotic powers of
specific tests to this envelope.

It is convenient to reparametrize our problem to $\theta=\sqrt{-\ln
(
1-h^{2}/c ) }$. As $h$ varies in the region of contiguity $ [
0,%
\sqrt{c} ) $, $\theta$ spans the entire half-line $ [ 0,\infty
) $. Note that the asymptotic mean and autocovariance functions
of the
log likelihood ratios derived in the previous section depend on $h$ only
through $h/\sqrt{c}=\sqrt{1-e^{-\theta^{2}}}$. Therefore, under the new
parametrization, they depend only on $\theta$. Loosely speaking,
$\theta$
and $\sqrt{p/n}\sim\sqrt{c}$ play the classical roles of a
``local parameter'' and a contiguity rate,
respectively.

Let $\beta( \theta_{1};\lambda) $ and $\beta( \theta
_{1};\mu) $ be the asymptotic powers of the asymptotically most
powerful $\lambda$- and $\mu$-based tests of size $\alpha$ of the
null $%
\theta=0$ against the alternative $\theta=\theta_{1}$. The following
proposition is proven in the \hyperref[app]{Appendix}.

\begin{proposition}\label{pr9}
Let $\Phi$ denote the standard normal distribution
function. Then
%
\begin{equation}
\label{localpower} \beta( \theta_{1};\lambda) =1-\Phi\biggl[
\Phi^{-1} ( 1-\alpha) -\frac{\theta_{1}}{\sqrt{2}} \biggr]
\end{equation}
and
%
\begin{equation}\label{localpowermu}
\beta( \theta_{1};\mu) = 1-\Phi\bigl[ \Phi^{-1} ( 1-
\alpha) -\sqrt{\tfrac{1}{2} \bigl( \theta_{1}^{2}-1+e^{-\theta
_{1}^{2}}
\bigr) } \bigr].
\end{equation}
\end{proposition}

\begin{figure}

\includegraphics{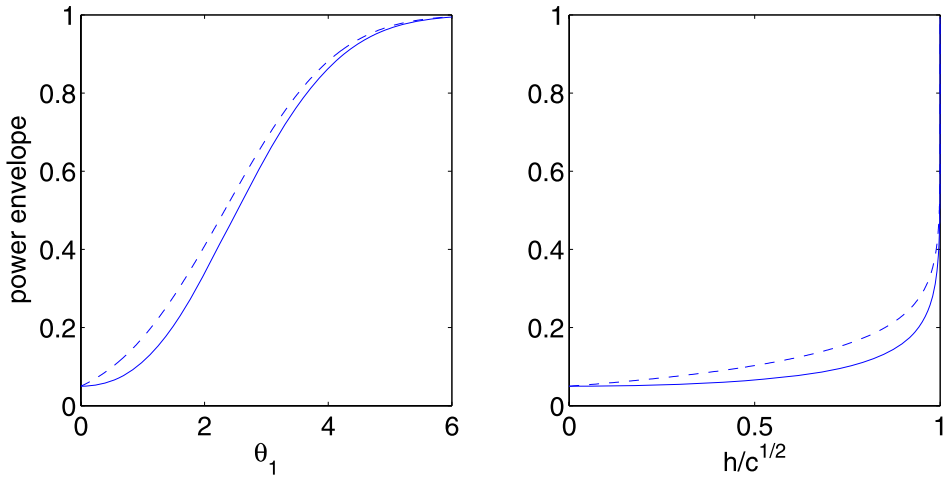}

\caption{The maximal asymptotic power of
the $\protect\lambda$-based tests (dashed lines) and $\protect\mu
$-based tests (solid lines)
of $\protect\theta=0$
against $\protect\theta=\protect\theta_{1}$. Left panel: $\protect
\theta$-parametrization. Right panel: $h$-parametrization.}\label{envelopes}
\end{figure}

Plots of the asymptotic power envelopes $\beta( \theta
_{1};\lambda
) $ and $\beta( \theta_{1};\mu) $ against $\theta_{1}$
for asymptotic size $\alpha=0.05$ are shown in the left panel of
Figure~\ref{envelopes}. The power loss of the $\mu$-based tests relative to the $%
\lambda$-based tests is due to the nonspecification of $\sigma^{2}$. In
contrast to $\lambda$-based tests, $\mu$-based tests may achieve the
corresponding power envelope even when $\sigma^{2}$ is unknown.\vadjust{\goodbreak}

The right panel of Figure~\ref{envelopes} shows the envelopes as functions of the original parameter $h$
normalized by $\sqrt{c}$. We see that the alternatives that can
theoretically be detected
with high probability are concentrated near the threshold $h=\sqrt{c}$.
The strong
nonlinearity of the $\theta$-parametrization should be kept in mind while
interpreting the figures that follow.

It is interesting to compare the power envelopes to the asymptotic
powers of
the likelihood ratio (LR) and weighted average power (WAP) tests. The $%
\lambda$-based LR and WAP tests of $\theta=0$ against the alternative
$%
\theta\in( 0,M ] $, where $M<\infty$, would reject the
null if
and only if, respectively, $2\sup_{\theta\in( 0,M ] }\ln
L (
\theta;\lambda) $ and $\ln\int_{0}^{M}L ( \theta;\lambda
) W ( \mathrm{d}\theta) $ are sufficiently large. The power
of a WAP test would, of course, depend on the choice of the weighting
measure $W( \mathrm{d}\theta) $. The $\mu$-based LR and WAP
tests are defined similarly. Theorem~\ref{th7} and Le Cam's third lemma
suggest a straightforward procedure for the numerical evaluation of the
corresponding
asymptotic power functions.

Consider, for example, the $\lambda$-based LR test statistic. According
to Theorem~\ref{th7}, its
asymptotic distribution under the null equals the distribution of
$2\sup_{\theta\in( 0,M ]} X_{\theta}$, where $X_{\theta}$ is
a Gaussian process with
$\mathrm{E} ( X_{\theta} ) =-\theta^{2}/4$ and $%
\operatorname{Cov} ( X_{\theta_{1}},X_{\theta_{2}} ) =-\frac
{1}{2}\ln
( 1-\sqrt{ ( 1-e^{-\theta_{1}^{2}} ) ( 1-e^{-\theta
_{2}^{2}} ) } ) $. According to Le Cam's third lemma, under a
specific alternative $\theta=\theta_{1}\leq M$, the asymptotic
distribution of the LR statistic equals the distribution of
$2\sup_{\theta\in( 0,M ] }\tilde{X}_{\theta}$, where $\tilde
{X}_{\theta}$
is a Gaussian process with the same covariance function as that of $%
X_{\theta}$, but with a different mean: $\mathrm{E} ( \tilde
{X}_{\theta
} ) =\mathrm{E} ( X_{\theta} ) + \operatorname{Cov}
( X_{\theta},X_{\theta_{1}} ) $.

Therefore, to numerically evaluate the asymptotic power function of the
$\lambda$-based LR test,
we simulate 500,000 observations of $X_{\theta}$ on a grid of 1000
equally spaced points in $\theta\in[ 0,M=6 ] $, where $M=6$
is chosen as the upper limit
of the grid because it is large enough for the power envelopes to rich
the value of 99\%.
For each observation, we save its supremum on the
grid, and use the empirical distribution of two times the suprema as the
approximate asymptotic distribution of the likelihood ratio statistic under
the null. We denote this distribution as $\hat{F}_{0}$.
Its 95\% quantile equals $4.3982$.

For each $\theta_{1}$ on the grid, we repeat the simulation
for process $\tilde{X}_{\theta}$ to obtain the approximate asymptotic
distribution of the likelihood ratio statistic under the alternative $%
\theta=\theta_{1}$, which we denote as $\hat{F}_{1}$. We use the
value of $\hat{F}_{1}$ at the 95\% quantile
of $\hat{F}_{0}$
as a numerical approximation
to the asymptotic power at $\theta_{1}$ of the $\lambda$-based LR test
with asymptotic size 0.05.

\begin{figure}

\includegraphics{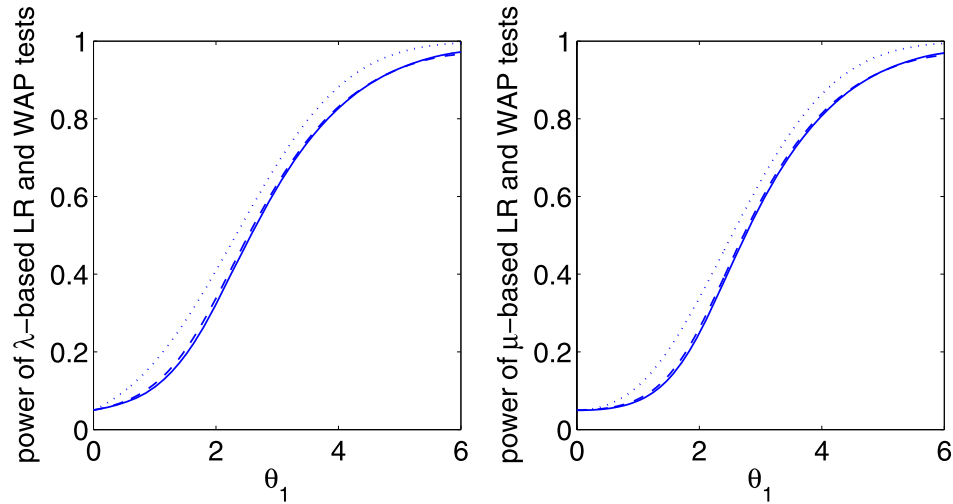}

\caption{The asymptotic power envelope
(dotted line), the asymptotic power of the LR test (solid line),
and the asymptotic power of the WAP test with uniform weighting
measure on $\protect\theta\in[ 0,6 ] $ (dashed line).
Left panel: $\protect\lambda$-based tests and envelope.
Right panel: $\protect\mu$-based tests and envelope.}\label{powerlr}
\end{figure}

Figure~\ref{powerlr} shows the resulting asymptotic power curve of the LR
test (solid line) along with the asymptotic power envelope (dotted
line). It
also shows the asymptotic power of the WAP test with $W ( \mathrm{d}
\theta) $ equal to the uniform measure on $ [ 0,6 ] $
(dashed line). The left and right panels correspond to $\lambda$- and
$\mu$-based
tests, respectively.

The asymptotic powers of the LR and WAP tests both come close to the
power envelope. The LR and WAP power functions are so close that they are
difficult to distinguish clearly. The asymptotic power of the WAP
test appears to be larger than that of the LR test for all $\theta
_{1}$ in
the $ [ 0,6 ] $ range, except for relatively large $%
\theta_{1}$. Hence, the LR test still may be admissible. More accurate
numerical analysis is needed to shed further light on this
issue.\looseness=-1

In the remaining part of this section, we consider some of the tests
that have been proposed previously in the literature, and, in
Proposition~\ref{pr10}, derive their asymptotic power functions. We focus on four
examples. Three of them are inspired by the
``classical'' fixed-$p$ theory, while the fourth is more
directly based on results from the large random matrix theory.

The problem of testing the hypothesis of sphericity has a long history,
and has generated a considerable body of literature, which we only very
briefly summarize here. The classical fixed-$p$ Gaussian analysis of
the various problems considered here goes back to \citet{Mau40},
who first derived the Gaussian likelihood ratio test for sphericity.
The (Gaussian) locally most powerful invariant (under shift, scale and
orthogonal transformations) test was obtained by John
\mbox{(\citeyear{Joh71}, \citeyear{Joh72})} and by \citet{Sug72}, with
adjusted versions resisting elliptical violations of the Gaussian
assumptions proposed in \citet{HalPai06}, where a Le Cam approach
is adopted under a general elliptical setting. \citet{LedWol02}
propose two extensions (for the unknown and known scale problems,
resp.) of John's test, while \citet{Baietal09} adapt Mauchly's
(\citeyear{Mau40}) likelihood ratio test.

\begin{example}[{[John's (\citeyear{Joh71}) test of sphericity]}]\label{ex1}
\citet{Joh71} proposes testing the
sphericity hypothesis $\theta=0$ against general alternatives using the
test statistic $U=\frac{1}{p}\operatorname{tr} [ ( \frac{\hat
{\Sigma}}{%
( 1/p ) \operatorname{tr} ( \hat{\Sigma} )
}-I_{p} ) ^{2}%
] $, where $\hat{\Sigma}$ is the sample covariance matrix of the data.
He shows that, when $n>p$, such a test is locally most powerful invariant.
Studying John's test when $p/n\rightarrow c\in(
0,\infty) $, \citet{LedWol02} prove that, under the null,
$nU-p\stackrel{d}{%
\rightarrow}N ( 1,4 ) $. Hence, the test with asymptotic size
$%
\alpha$ rejects the null hypothesis of sphericity if $\frac{1}{2}
( nU-p-1 )
>\Phi^{-1}( 1-\alpha) $.
\end{example}

\begin{example}[{[The \citet{LedWol02} test of $\Sigma=I$]}]\label{ex2}
\citet{LedWol02}
propose using $W=\frac{1}{p}\operatorname{tr} [ ( \hat{\Sigma
}-I )
^{2} ] -\frac{p}{n} [ \frac{1}{p}\operatorname{tr}\hat{\Sigma} ]
^{2}+\frac{p%
}{n}$ as a test statistic for testing the hypothesis that the population
covariance matrix is a unit matrix. Under the null, $nW-p\stackrel{d}{%
\rightarrow}N ( 1,4 ) $. As in the previous example, the null
hypothesis is
rejected at asymptotic size $\alpha$ if $\frac{1}{2}%
( nW-p-1 ) >\Phi^{-1}( 1-\alpha) $.
\end{example}

\begin{example}[{[The ``corrected'' LRT of \citet{Baietal09}]}]\label{ex3}
When $n>p$, \citet{Baietal09} propose\vspace*{1pt} a corrected version of
the likelihood ratio statistic $\mathrm{CLR}=\operatorname{tr}\hat{\Sigma}-\ln
\det\hat{%
\Sigma}-p-p ( 1- ( 1-\frac{n}{p} ) \ln( 1-\frac{p}{n}
) ) $ based on the entire data, as opposed to $\lambda$ or
$\mu
$ only, to test the equality of the population covariance matrix to the
identity matrix against general alternatives. Under the null,
$\mathrm{CLR}\stackrel{d}{%
\rightarrow}N ( -\frac{1}{2}\ln( 1-c ),-2\ln(
1-c ) -2c ) $. The null hypothesis is rejected whenever $\mathrm{CLR}
+ \frac{1}{2}\ln( 1-c ) >
\sqrt{-2\ln(
1-c ) -2c} \Phi^{-1} ( 1-\alpha)$.\vadjust{\goodbreak}
\end{example}

More directly inspired by the asymptotic theory of random matrices,
several authors have recently proposed and studied various tests based
on $\lambda _{1}$ or $\mu_{1}$: see \citet{Bej05},
\citet{PatPriRei06}, \citet{KriNad09},
\citet{Ona09}, \citet{Biaetal} and \citet{NadSil10}. We
refer to these tests, which reject $H_{0}$ for large values of
$\lambda_{1}$ or $\mu_{1}$, as Tracy--Widom-type tests.

\begin{example}[(Tracy--Widom-type tests)]\label{ex4}
Asymptotic critical values of such tests
are obtained using the fact, established by \citet{Joh01}, that
under the
null,
%
\begin{equation}\label{TW}
{n^{{2}/{3}}c^{{1}/{6}}} { ( 1+\sqrt{c} )
^{-{4}/{3}}%
} \bigl( \lambda_{1}- ( 1+\sqrt{c} )
^{2} \bigr) \stackrel{d} {%
\rightarrow}\mathrm{TW},
\end{equation}
where TW denotes the Tracy--Widom law of the first kind. The null hypothesis
is rejected when $\lambda_{1}$ or $\mu_{1}$ exceeds the adequate
Tracy--Widom quantile.
\end{example}

Consider the tests described in Examples~\ref{ex1},~\ref{ex2},
\ref{ex3} and~\ref{ex4}, and denote by $\beta_{\mathrm{J}} ( \theta_{1} ) $,
$\beta_{\mathrm{LW}} ( \theta_{1} ) $, $\beta_{\mathrm{CLR}} ( \theta_{1} ) $, and
$\beta_{\mathrm{TW}} ( \theta_{1} ) $ their respective asymptotic powers at
asymptotic level $\alpha$. The following proposition is established in
the \hyperref[app]{Appendix}.

\begin{proposition}\label{pr10}
Denote $1-e^{-\theta_{1}^{2}}$ as $\psi( \theta
_{1} ) $. The asymptotic power functions of the tests described
in Examples~\ref{ex1}--\ref{ex4} satisfy, for any $\theta_1 >0$,
%
\begin{eqnarray}
\label{powerTW} \beta_{\mathrm{TW}} ( \theta_{1} ) &=&\alpha,
\\
\label{betaJ}
\beta_{\mathrm{J}} ( \theta_{1} ) &=&\beta_{\mathrm{LW}} (
\theta_{1} ) =1-\Phi\bigl( \Phi^{-1} ( 1-\alpha) -
\tfrac{1}{2}\psi( \theta_{1} ) \bigr)
\end{eqnarray}
and
%
\begin{equation}\label{betaCCLR}\quad
\beta_{\mathrm{CLR}} ( \theta_{1} ) =1-\Phi\biggl(
\Phi^{-1} ( 1-\alpha) -\frac{\sqrt{c\psi( \theta_{1} ) }-\ln
( 1+%
\sqrt{c\psi( \theta_{1} ) } ) }{\sqrt{-2\ln(
1-c ) -2c}} \biggr).
\end{equation}
\end{proposition}

With the important exception of \citet{Sri05}, (\ref
{powerTW})--(\ref{betaCCLR}) are the first
results on the asymptotic power of those tests against contiguous
alternatives. \citet{Sri05} analyzes the asymptotic power of tests
similar to those in Examples~\ref{ex1} and~\ref{ex2}. His Theorems 3.1 and 4.1 can be
used to establish (\ref{betaJ}).

From Proposition~\ref{pr10}, we see that the local asymptotic power of the
Tracy--Widom-type tests is trivial. As shown by \citet{BaiBenPec05} in the
complex data case and by \citet{FerPec09} in the real data
case, the convergence (\ref{TW}) holds not only under the null, but also
under any alternative of the form $h=h_{0}<\sqrt{c}$. Under the
``local'' parametrization adopted in this
section, such alternatives have the form $\theta=\theta_{1}>0$. It
can be
shown that the Tracy--Widom-type tests are consistent against noncontiguous
alternatives $h=h_1 >\sqrt{c}$. However, such a consistency is likely
to be
also a property of\vadjust{\goodbreak} the LR tests based on $\mu$ or on $\lambda$. If this
holds true, the LR tests asymptotically dominate the Tracy--Widom-type tests.
A more detailed analysis of the optimality properties of LR
tests is the subject of ongoing research.

\begin{figure}

\includegraphics{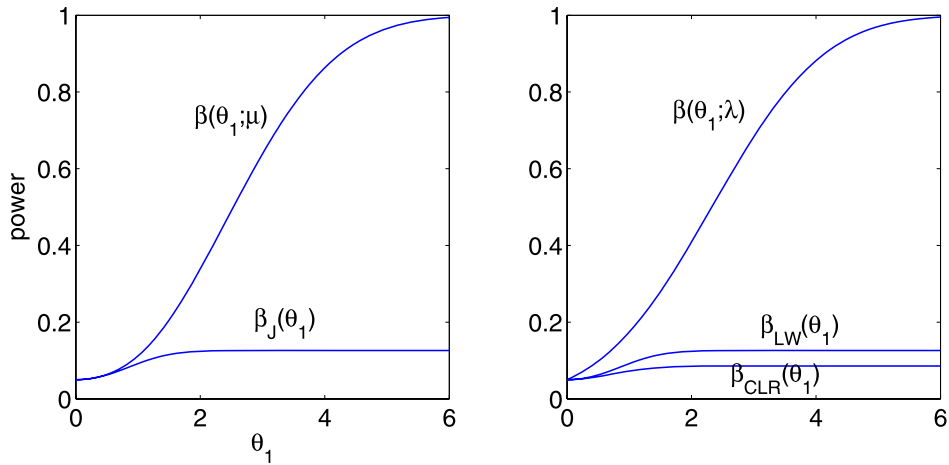}

\caption{Asymptotic powers $ ( \beta_{\mathrm{J}}, \beta_{\mathrm{LW}},
\beta _{\mathrm{CLR}} )$ of the tests described in Examples \protect\ref{ex1}
(John), \protect\ref{ex2} (Ledoit and Wolf) and~\protect\ref{ex3} (Bai et
al.).}\label{powerspecificR}
\end{figure}

The asymptotic power functions of the tests from Examples~\ref{ex1},
\ref{ex2} and~\ref{ex3} are nontrivial. Figure~\ref{powerspecificR}
compares these power functions to the corresponding power envelopes.
Since John's test is invariant with respect to orthogonal
transformations and scalings, $\beta_{\mathrm{J}} ( \theta _{1} ) $ is compared
to the power envelope $\beta( \theta_{1};\mu)$. The asymptotic power
functions $\beta_{\mathrm{LW}} ( \theta_{1} ) $
and $%
\beta_{\mathrm{CLR}} ( \theta_{1} ) $ are compared to the power
envelope $%
\beta( \theta_{1};\lambda) $ because the Ledoit--Wolf
test of $%
\Sigma=I$ and the ``corrected'' likelihood
ratio test are invariant only with respect to orthogonal transformations.

Interestingly, whereas $\beta_{\mathrm{J}} ( \theta_{1} ) $ and $\beta
_{\mathrm{LW}} ( \theta_{1} ) $ depend only on $\alpha$ and $\theta_{1}$,
$\beta_{\mathrm{CLR}} ( \theta_{1} ) $ depends also on $c$. As $c$
converges to one, $\beta_{\mathrm{CLR}} ( \theta_{1} ) $ converges
to $\alpha$, which corresponds to the case of trivial power. As $c$ converges
to zero, $\beta_{\mathrm{CLR}} ( \theta_{1} ) $ converges to $\beta
_{\mathrm{J}} ( \theta_{1} ) $. In Figure~\ref{powerspecificR}, we provide
the plot of $\beta_{\mathrm{CLR}} ( \theta_{1} ) $ that corresponds
to $%
c=0.5$.

The left panel of Figure~\ref{powerspecificR} shows that the power function
of John's test is very close to the power envelope $\beta( \theta
_{1};\mu) $ in the vicinity of $\theta_{1}=0$. Such behavior is
consistent with the fact that John's test is locally most powerful
invariant. However, for large $\theta_{1}$, the asymptotic power functions
of all the tests from Examples~\ref{ex1},~\ref{ex2} and~\ref{ex3} are lower than the corresponding
asymptotic power envelopes. We should stress here that these tests have
power against general alternatives as opposed to the
``spiked'' alternatives that maintain the assumption
that the
population covariance matrix of data has the form $\sigma^{2} (
I_{p}+hvv^{\prime} ) $.\vadjust{\goodbreak}

For the ``spiked'' alternatives, the $\lambda$- and $\mu$-based
LR tests
may be more attractive. However, implementing these tests requires
some care. A ``quick-and-dirty'' approach would be to approximate
$\ln L (\theta;\lambda)$ and $\ln L (\theta;\mu
)$ by
the simple but asymptotically equivalent expressions from (\ref
{equivalence1})
and (\ref{equivalence2}), compute two times their maxima on a grid
over $\theta\in( 0,M ]$, and compare them with critical
values obtained
by simulation as for the construction of Figure~\ref{powerlr}.
Unfortunately, in finite samples, this simple approach will lead to a
numerical breakdown
whenever $z_{0}(h(\theta))$ happens to be less than the largest sample
covariance eigenvalue
for some $\theta\leq M$. In addition, since the asymptotic
approximation derived in
Theorem~\ref{th7} is not uniform over entire half-line $\theta\in
[0,\infty)$,
its quality will depend on the choice of $M$. For relatively large $M$,
the asymptotic
behavior of the LR test implemented as above may poorly match its
finite sample behavior.

Instead, we recommend implementing the LR tests without using the
asymptotic approximations.
The finite sample log likelihood ratios $\ln L (\theta;\lambda
)$ and
$\ln L (\theta;\mu)$ can be computed using the contour
integral representations
(\ref{LRcontourr1}) and (\ref{LRcontourr2}).
Choosing the contour of integration so that the sample covariance
eigenvalues remain
to its left will eliminate the numerical breakdown problem associated
with the asymptotic
tests. Furthermore, under the Gaussianity assumption,
the finite sample distributions of the log likelihood ratios are
pivotal. Hence, the exact
critical values can be computed via Monte Carlo simulations as follows:
simulate many replications
of data under the null. For each replication, compute the log
likelihood ratio and store two times
its maximum.
Use the 95\% quantile of the empirical distribution of the stored
values as a numerical
approximation for the exact critical value of the test. The finite
sample properties
of such a test are left as an important topic for future research.

\section{Conclusion}\label{sec6}

In this paper, we study the asymptotic power of tests for the existence of
rank-one perturbations of sphericity as both the dimensionality of the data
and the number of observations go to infinity. Focusing on tests that are
invariant with respect to orthogonal transformations and rescaling, we
establish the
convergence of the log ratio of the joint densities of the sample covariance
eigenvalues under the alternative and null hypotheses to a Gaussian process
indexed by the norm of the perturbation.

When the perturbation norm is larger than the phase transition threshold
studied in \citet{BaiBenPec05}, the limiting log-likelihood process is
degenerate and the joint eigenvalue distributions under the null and
alternative hypotheses are asymptotically mutually singular, so that the
discrimination between the null and the alternative is asymptotically
certain. When the norm is below the threshold, the limiting log-likelihood
process is nondegenerate and the joint eigenvalue distributions under the
null and alternative hypotheses are mutually contiguous. Using the asymptotic
theory of statistical experiments,\vadjust{\goodbreak} we obtain power envelopes and derive
the asymptotic
size and power for various eigenvalue-based tests in the region of
contiguity.

Several questions are left for future research. First, we only considered
rank-one perturbations of the spherical covariance matrices. It would
be desirable to extend the analysis to finite-rank perturbations. Such an
extension will require a more complicated technical analysis. Second, it
would be interesting to extend our analysis to the asymptotic
regime $p,n\rightarrow\infty$ with $p/n\rightarrow\infty$ or
$p/n\rightarrow0$. In the context of sphericity tests, such asymptotic
regimes have been recently studied in \citet{BirDet05}.
Third, a thorough analysis of the finite sample properties of the
proposed LR tests would clarify the
related practical implementation issues.
Fourth, our Lemma~\ref{le5} can be used to derive
higher-order asymptotic approximations to the likelihood ratios, which may
improve finite-sample performances of asymptotic tests. Finally, it
would be of
considerable interest to relax the Gaussian assumptions, for example, into
elliptical ones, preferably with unspecified radial densities, on the model
(in a fixed-$p$ context) of \citet{HalPai06}.\looseness=-1

\begin{appendix}\label{app}
\section*{Appendix}
\subsection{\texorpdfstring{Proof of Proposition \protect\ref{pr1}}{Proof of Proposition 1}}\label{sec7.1}

For the joint density $p ( \lambda;h ) $ of $\lambda
_{1},\ldots,\lambda_{m}$, we have%
%
\begin{eqnarray}\label{commonreal}
p ( \lambda;h ) &=&\tilde{\gamma}\frac{\prod_{i=1}^{m}\lambda
_{i}^{%
({\llvert p-n\rrvert-1})/{2}}\prod_{i<j}^{m} ( \lambda
_{i}-\lambda_{j} ) }{ ( 1+h ) ^{n/2}}\nonumber\\[-8pt]\\[-8pt]
&&{}\times\int_{\mathcal{O}%
(p)}e^{-({n}/{2})\operatorname{tr} ( \Pi Q^{\prime}\Lambda
Q )}( \mathrm{d}Q ),\nonumber
\end{eqnarray}
where $\tilde{\gamma}$ depends only on $n$ and $p$, $\Pi=\operatorname{diag}
( ( 1+h ) ^{-1},1,\ldots,1 ) $, $\mathcal{O}(p)$ is
the set
of all $p\times p$ orthogonal matrices and $ ( \mathrm{d}Q ) $ is
the invariant measure on the orthogonal group $\mathcal{O}(p)$
normalized to
make the total measure unity. When $n\geq p$, (\ref{commonreal}) is
a special case of the density given in \citet{Jam64}, page 483. When $%
n<p$, (\ref{commonreal}) follows from Theorems 2 and 6 in \citet{Uhl94}.

Let $\Psi=\operatorname{diag} ( \frac{h}{1+h},0,\ldots,0 ) $ be
a $p\times
p$ matrix. Since $\Pi=I_{p}-\Psi$, we have $\operatorname{tr} (
\Pi Q^{\prime
}\Lambda Q ) =\operatorname{tr}\Lambda-\operatorname{tr} (
\Psi Q^{\prime
}\Lambda Q )$, and we can rewrite (\ref{commonreal}) as%
%
\begin{eqnarray}\label{LRreal1}
p ( \lambda;h ) &=&\tilde{\gamma}\frac{\prod_{i=1}^{m}\lambda
_{i}^{({\llvert p-n\rrvert-1})/{2}} \prod_{i<j}^{m} ( \lambda
_{i} - \lambda_{j} ) e^{-({n}/{2})\operatorname{tr}\Lambda
}}{ (
1+h ) ^{n/2}} \nonumber\\[-8pt]\\[-8pt]
&&{}\times\int_{\mathcal{O}(p)}e^{({n}/{2})\operatorname{tr}
( \Psi Q^{\prime}\Lambda Q ) }
( \mathrm{d}Q ).\nonumber
\end{eqnarray}
Note that $\operatorname{tr} ( \Psi Q^{\prime}\Lambda Q )
=\operatorname{tr}
( Q\Psi Q^{\prime}\Lambda) =\frac{h}{1+h}x_{p}^{\prime
}\Lambda
x_{p}$, where $x_{p}$ is the first column of $Q$. When $Q$ is uniformly
distributed over $\mathcal{O}(p)$, its first column $x_{p}$ is uniformly
distributed over $\mathcal{S} ( p ) $. Therefore, we have%
%
\begin{eqnarray}\label{detailedp2}
p ( \lambda;h ) &=&\tilde{\gamma}\frac{\prod_{i=1}^{m}\lambda
_{i}^{({\llvert p-n\rrvert-1})/{2}} \prod_{i<j}^{m} ( \lambda
_{i} - \lambda_{j} ) e^{-({n}/{2})\operatorname{tr}\Lambda
}}{ (
1+h ) ^{n/2}} \nonumber\\[-8pt]\\[-8pt]
&&{}\times\int_{\mathcal{S} ( p ) }e^{({n}/{2})
({h}/({1+h}))x_{p}^{\prime}\Lambda x_{p}}
( \mathrm{d}x_{p} ),\nonumber
\end{eqnarray}
which establishes (\ref{alternativereal}). Now, let $y=\lambda
_{1}+\cdots+\lambda_{p}$ so that $\mu_{j}=\lambda_{j}/y$. Note that $%
\operatorname{tr}\Lambda=y$, $\operatorname{tr}M= \mu_{1}+\cdots+\mu
_{p}=1$, and
that the Jacobian of the coordinate change from $\lambda
_{1},\ldots,\lambda_{m}$
to $\mu_{1},\ldots,\mu_{m-1},y$ equals $y^{m-1}$. Changing variables in
(\ref{detailedp2}), and integrating $y$ out, we obtain (\ref{alternativerealmu}).

\subsection{\texorpdfstring{Proof of Lemma \protect\ref{le3}}{Proof of Lemma 3}}\label{sec7.2}

Using (\ref{Erdelyirepresentation}) in the ratio of the right-hand
side of (%
\ref{alternativereal}) with $h>0$ to that with $h=0$, and changing the
variable of integration from $s$ to $z=\frac{1+h}{h}\frac{2}{n}s$, we
get (%
\ref{LRcontourr1}). Further, from (\ref{alternativerealmu}), we
have%
%
\begin{eqnarray}\label{mudensityh0}
p ( \mu;0 ) &=&\delta( n,p,\mu) \int_{0}^{\infty
}y^{{np}/{2}-1}e^{-{ny}/{2}}
\,\mathrm{d}y\nonumber\\[-8pt]\\[-8pt]
&=&\delta( n,p,\mu) \biggl( \frac{2}{n} \biggr) ^{{np}/{2}}
\Gamma\biggl( \frac
{np}{2} \biggr).\nonumber
\end{eqnarray}
For $h>0$, using (\ref{Erdelyirepresentation}) in (\ref%
{alternativerealmu}), we get%
\begin{eqnarray*}
p ( \mu;h ) &=&\frac{\delta( n,p,\mu) }{ (
1+h ) ^{n/2}}\frac{\Gamma( p/2 ) }{2\pi i}\\
&&{}\times\int_{0}^{\infty
}
\oint_{\tilde{\mathcal{K}}}y^{({np-2})/{2}}e^{s-{n}y/{2}}\prod
_{j=1}^{p} \biggl( s-\frac{n}{2}
\frac{yh}{1+h}\mu_{j} \biggr) ^{-%
{1}/{2}}\,\mathrm{d}s
\,\mathrm{d}y,
\end{eqnarray*}
where $\tilde{\mathcal{K}}$ is a contour starting at $-\infty$, encircling
counter-clockwise the points~$0$, $\frac{ny}{2}\frac{h}{1+h}\mu
_{1},\ldots,
\frac{ny}{2}\frac{h}{1+h}\mu_{m}$ and going back to $-\infty$. Since
$%
\frac{h}{1+h}\mu_{j}<1$ by construction, we may and will choose
$\tilde{\mathcal{%
K}}$ so that for any $s\in\tilde{\mathcal{K}}$, $\operatorname
{Re}s<\frac{ny%
}{2}$. Changing variables of integration from $y$ and $s$ to $w=\frac{ny}{2}$
and $z=s\frac{1+h}{hw}S$, where $S$ is any positive constant, and dividing
by the right-hand side of (\ref{mudensityh0}), we obtain%
\begin{eqnarray*}
L ( h;\mu) &=& \frac{S^{({p-2})/{2}} ( 1+h ) ^{({p-n-2})/{2}}\Gamma(
{p}/{2} ) }{h^{({p-2})/{2}}\Gamma(
{np}/{2} ) 2\pi i}\\
&&{}\times\int_{0}^{\infty}
\oint_{\mathcal{K}}w^{{np}/{2}-{p}/{2}}e^{({wh}/({1+h}))({z}/{S})-w}\prod%
_{j=1}^{p} ( z-S\mu_{j} ) ^{-{1}/{2}}
\,\mathrm{d}z\,\mathrm{d}w,
\end{eqnarray*}
where $\mathcal{K}$ is a contour starting at $-\infty$, encircling
counter-clockwise the points~$0$, $S\mu_{1},\ldots,S\mu_{m}$, and going
back to $%
-\infty$. In addition, for any $z\in\mathcal{K}$, $\operatorname
{Re}z<\frac{1+h}{h}%
S$. Such a choice of $\mathcal{K}$ guarantees that the integrand in the
above double integral is absolutely integrable on $ [ 0,\infty
)
\times\mathcal{K}$, so that Fubini's theorem can be used to justify the
interchange of the order of the integrals. Changing the order of the
integrals and setting $S=\lambda_{1}+\cdots+\lambda_{p}$, we obtain (\ref
{LRcontourr2}).

\subsection{\texorpdfstring{Proof of Theorem \protect\ref{th7}}{Proof of Theorem 7}}\label{sec7.3}

First, let us formulate the following technical lemma. Its proof is in the
Supplementary Appendix [\citet{OnaMorHal}].

\begin{lemma}\label{le11}
\textup{(i)} If $h<\sqrt{c_{p}}$, ${f_{0}=-\frac{1}{2} (
c_{p}+ (
1-c_{p} ) \ln( 1+h ) -c_{p}\ln\frac{c_{p}}{h} ).}$

\textup{(ii)} If $h>\sqrt{c_{p}}$, ${f_{0}=-\frac{1}{2} (
h+c_{p}+ (
1-c_{p} ) \ln( c_{p}+h ) -\frac{c_{p}}{h}-\ln h ).}$

\end{lemma}

Below, we prove Theorem~\ref{th7} for $L ( h;\mu) $. The proof for
$L (
h;\lambda) $ is similar but simpler, and we omit it to save
space. As
follows from Lemmas~\ref{le4} and~\ref{le5}, the integral in (\ref{LRcontourr2}) can be
represented as $2e^{-nf_{0}} [ \Gamma( \frac{1}{2} )
\frac{%
a_{0}}{n^{1/2}}+\frac{O_{p} ( 1 ) }{hn^{3/2}} ] $
uniformly in
$h\in( 0,\bar{h} ] $. Therefore, and since $\Gamma(
\frac{1%
}{2} ) =\sqrt{\pi}$, we can write%
%
\begin{equation}\label{likelihoodmuincontiguity}
L ( h;\mu) =\frac{k_{2}S^{({p-2})/{2}}}{\sqrt{n\pi}i}%
e^{-nf_{0}} \biggl[
a_{0}+h^{-1}O_{p} \biggl( \frac{1}{n}
\biggr) \biggr],
\end{equation}
where $k_{2}=h^{-({p-2})/{2}} ( 1+h ) ^{({p-n-2})/{2}}\frac
{%
( n-1 ) p}{2}\Gamma( \frac{ ( n-1 ) p}{2} )
\Gamma( \frac{p}{2} ) \Gamma^{-1} ( \frac{np}{2} ) $.
Using Stirling's approximation $\Gamma( r ) =e^{-r}r^{r} (
\frac{2\pi}{r} ) ^{1/2} ( 1+O ( r^{-1} ) ) $
with $%
r=\frac{p}{2}$, $\frac{np}{2}$ and $\frac{ ( n-1 ) p}{2}$, and the
fact that $\ln( n-1 ) =\ln n-n^{-1}-\frac{1}{2}n^{-2}+O (
n^{-3} ) $, we find, after algebraic simplifications, that%
%
\begin{eqnarray}\label{gammas}
\frac{k_{2}}{\sqrt{n\pi}}&=&h^{-({p-2})/{2}} ( 1+h )
^{({p-n-2})/{2}}\nonumber\\[-8pt]\\[-8pt]
&&{}\times e^{-(({p-2})/{2})\ln
n-{p}/{2}+{c_{p}}/{4}+\ln c_{p}/2} \bigl( 1+O \bigl(
n^{-1} \bigr) \bigr).\nonumber
\end{eqnarray}

Using (\ref{gammas}) and Lemma~\ref{le11}(i), we obtain
\[
\frac{k_{2}S^{({p-2})/{2}}}{\sqrt{n\pi}i}e^{-nf_{0}}h^{-1}O_{p} \biggl(
\frac{1}{n} \biggr) =\frac{1}{1+h} \biggl( \frac{S}{p} \biggr)
^{({p-2})/{2}%
}e^{{c_{p}}/{4}-\ln c_{p}/2}O_{p} \biggl( \frac{1}{n}
\biggr),
\]
which, together with the fact that $S-p=O_{p} ( 1 ) $, implies
that%
%
\begin{equation}\label{theor1part1}
\frac{k_{2}S^{({p-2})/{2}}}{\sqrt{n\pi}i}e^{-nf_{0}}h^{-1}O_{p} \biggl(
\frac{1}{n} \biggr) =O_{p} \biggl( \frac{1}{n} \biggr)
\end{equation}
uniformly over $h\in( 0,\bar{h} ] $.

Now, as can be verified using (\ref{fdefinition}) and (\ref
{Stijeltjesanalytic}), if $h<\sqrt{c_{p}}$, then%
%
\begin{equation}\label{f2}
f_{2}=-\frac{h^{2}}{4 ( 1+h ) ^{2} ( c_{p}-h^{2} ) }.
\end{equation}
Therefore, using (\ref{twocoefficients}), we obtain%
%
\begin{equation}\label{a0}
a_{0}=i\frac{ ( 1+h ) ( c_{p}-h^{2} ) ^{1/2}}{h}g_{0}.
\end{equation}
Using (\ref{gdefinition}), (\ref{gammas}), (\ref{a0}) and Lemma~\ref{le11}(i)
in (%
\ref{likelihoodmuincontiguity}), after algebraic simplifications and
rearrangements of terms, we get%
%
\begin{eqnarray}\label{longone}\quad
&&
\ln\biggl[ \frac{k_{2}S^{
({p-2})/{2}}e^{-nf_{0}}a_{0}}{\sqrt{n\pi}i}%
\biggr]\nonumber \\
&&\qquad=
\frac{1}{2}\ln\biggl( 1-\frac{h^{2}}{c_{p}} \biggr) +\frac
{c_{p}}{4}+%
\frac{p-2}{2}\ln\biggl( \frac{S}{p} \biggr)
\\
&&\qquad\quad{}-\frac{n}{2}\frac{hz_{0} ( h ) }{1+h}-\frac
{np-p+2}{2}\ln\biggl(
1-%
\frac{h}{1+h}\frac{z_{0} ( h ) }{S} \biggr) -\frac{1}{2}
\Delta_{p} \bigl( z_{0} ( h ) \bigr).\nonumber
\end{eqnarray}
Finally, using the fact that $S-p=O_{p} ( 1 ) $, we obtain $\ln
( {S}/{p} ) =
{(S-p)}/{p}+O_{p} ( p^{-2} )$ and
\begin{eqnarray*}
\ln\biggl( 1 - \frac{h}{1 + h}\frac{z_{0} (h)}{S} \biggr) &=& -%
\frac{h}{1 + h}\frac{z_{0} (h)}{p}-\frac{1}{2} \biggl( \frac{hz_{0} (h)}{
( 1 + h )p}
\biggr) ^{2}\\
&&{} + \frac{h}{ 1 + h }\frac{z_{0} (h)}{p^{2}} ( S - p )
+ O_{p} \bigl( p^{-3} \bigr).
\end{eqnarray*}

The latter two equalities, (\ref{longone})
and the fact that $\frac{h}{1+h}z_{0} ( h ) =h+c_{p}$ entail%
%
\begin{eqnarray}\label{theor1part2}
&&
\frac{k_{2}S^{({p-2})/{2}}e^{-nf_{0}}a_{0}}{\sqrt{n\pi}i}\nonumber\\[-8pt]\\[-8pt]
&&\qquad=e^{-
\{ \Delta_{p} ( z_{0} ( h ) ) -\ln(
1-{h^{2}}/{c_{p}} ) +({h}/{c_{p}}) ( S-p ) -
{h^{2}}/({2c_{p}})%
+O_{p} ( p^{-1} ) \}/2 },\nonumber
\end{eqnarray}
which, together with (\ref{theor1part1}), imply formula (\ref{equivalence2}).

Now, let us prove the convergence of $\ln L ( h;\mu) $ to $%
\mathcal{L} ( h;\mu) $. By (\ref{equivalence2}), the joint
convergence of $\ln L ( h_{j};\mu) $ with $j=1,\ldots,r$ to a
Gaussian vector is equivalent to the convergence of $ ( S-p,\Delta
_{p} ( z_{0}(h_{1}) ),\ldots,\Delta_{p} ( z_{0}(h_{r}) )
) $ to a Gaussian vector. A proof of the following technical lemma,
based on Theorem 1.1 of \citet{BaiSil04}, is given in the
Supplementary Appendix [\citet{OnaMorHal}].

\begin{lemma}\label{le12}
Suppose that the null hypothesis holds. Then, as $p,n\rightarrow_{c}
\infty$, the vector $ ( S-p,\Delta_{p} ( z_{0}(h_{1}) ),\ldots,\Delta
_{p} ( z_{0}(h_{r}) ) ) $ converges in distribution to a
Gaussian vector $ ( \eta,\xi_{1},\ldots,\xi_{r} ) $ with
\begin{eqnarray*}
\mathrm{E}\eta&=&0,\qquad \operatorname{Var}(\eta) =2c,\qquad
\operatorname{Cov}(\eta,\xi_{j} ) =-2h_{j},
\\
\operatorname{Cov}(\xi_{j},\xi_{k}) &=&-2\ln\bigl(
1-c^{-1}h_{j}h_{k} \bigr)\quad \mbox{and}\quad \mathrm{E}
\xi_{j}=\tfrac{1}{2}\ln\bigl( 1-c^{-1}h_{j}^{2}
\bigr).
\end{eqnarray*}
\end{lemma}

Lemma~\ref{le12} and (\ref{equivalence2}) imply that
$\mathrm{E} [ \mathcal{L} ( h_{j};\mu) ] =-\frac
{1}{2}%
\mathrm{E}\xi_{j}+\frac{1}{2}\ln( 1-c^{-1}h_{j}^{2} ) +\frac
{1}{4%
}c^{-1}h_{j}^{2}=\frac{1}{4} [ \ln( 1-c^{-1}h_{j}^{2} )
+c^{-1}h_{j}^{2} ]$ and%
\begin{eqnarray*}
\operatorname{Cov} \bigl[ \mathcal{L} ( h_{j};\mu),\mathcal{L} (
h_{k};\mu) \bigr] &=& \frac{1}{4}\operatorname{Cov} ( \xi
_{j},\xi_{k} ) +\frac{h_{k}}{4c}\operatorname{Cov} (
\xi_{j},\eta)\\
&&{} +
\frac{h_{j}}{4c}\operatorname{Cov}(\xi_{k},\eta)
+\frac{h_{j}h_{k}}{4c^{2}}\operatorname{Var} ( \eta) \\
&=&-\frac{1}{2}%
\ln\bigl( 1 - c^{-1}h_{j}h_{k} \bigr) -
\frac{h_{j}h_{k}}{2c},
\end{eqnarray*}
which establishes (\ref{meanmu}) and (\ref{covariancemu}).

To complete the proof of Theorem~\ref{th7}, we need to note that the tightness
of $L ( h;\mu) $, viewed as a random element of the space $C
( %
[ 0,\bar{h} ] )$, as $p,n\rightarrow_{c} \infty$,
follows from formula (\ref{equivalence2}) and the fact
that $S-p
$ and $\Delta_{p} ( z_{0} ( h )
) 
$, are $O_{p}(1)$, uniformly in $h\in( 0,\bar{h}%
] $. This uniformity is a consequence of Lemma A2 proven in the
Supplementary Appendix [\citet{OnaMorHal}].

\subsection{\texorpdfstring{Proof of Theorem \protect\ref{th8}}{Proof of Theorem 8}}\label{sec7.4}

As in the proof of Theorem~\ref{th7}, we will focus on the case of the likelihood
ratio based on $\mu$. The proof for $L ( h;\lambda) $ is similar.
According to Lemma~\ref{le6} and formula (\ref{LRcontourr2}), for any $\tilde
{h}>%
\sqrt{c}$, we have $L ( h;\mu) =k_{2}S^{({p-2})/{2}%
}e^{-nf ( z_{0} ( \tilde{h} ) ) }O_{p} ( 1
) $.
Using (\ref{gammas}) and the fact that $ ( \frac{S}{p} )
^{p}= ( 1+\frac{S-p}{p} ) ^{p}= ( 1+\frac{O_{p} (
1 )
}{p} ) ^{p}=O_{p} ( 1 ) $, we can write
%
\begin{equation}\label{firstsimplification}
L ( h;\mu) =e^{({n}/{2}) ( c_{p}\ln({c_{p} (
1+h ) }/{h})-\ln( 1+h ) -c_{p}-2f ( z_{0} ( \tilde
{h}%
) ) ) }O_{p} \bigl( n^{1/2} \bigr).
\end{equation}
Noting that $\tilde{h}>\sqrt{c_{p}}$ for sufficiently large $n$ and
$p$, and
using Lemma~\ref{le11}(ii) and the fact that $\frac{\tilde{h}}{1+\tilde{h}}%
z_{0} ( \tilde{h} ) =\tilde{h}+c_{p}$, we get
$-2f ( z_{0} ( \tilde{h} ) ) = ( 1 - c_{p}
) \ln( c_{p}+%
\tilde{h} ) -\frac{c_{p}}{\tilde{h}}-\ln\tilde{h}+\frac{h}{1+h}%
z_{0} ( \tilde{h} )$.
Substituting the latter expression in (\ref{firstsimplification}) and simplifying, we obtain%
%
\begin{equation}\label{LRasymptotica}
L ( h;\mu) =e^{({n}/{2})R ( h,\tilde{h},c_{p} )
}O_{p} \bigl( n^{1/2} \bigr),
\end{equation}
where $O_{p} ( \cdot) $ is uniform in $h\in[ \tilde{h},\infty) $ and $R
( h,\tilde{h},c_{p} ) = ( 1 - c_{p} ) \ln( c_{p} + \tilde{h} )
-\frac{c_{p}}{\tilde{h}}-\ln\tilde{h}+\frac
{h}{1+h}%
z_{0} ( \tilde{h} )- ( 1-c_{p} ) \ln( 1+h
) -c_{p}\ln h+c_{p}\ln
c_{p}-c_{p}$.

As $n,p\rightarrow\infty$, $R ( h,\tilde{h},c_{p} )
\rightarrow
R ( h,\tilde{h},c ) $ uniformly over $ ( h,\tilde{h}
) \in%
[ \sqrt{c},H ] ^{2}$. On the other hand, $R ( h,\tilde{h},c ) $
is continuous on $ ( h,\tilde{h} ) \in[ \sqrt{c},H%
] ^{2}$,
$R ( \sqrt{c},\sqrt{c},c ) = 0$, and $\frac{\mathrm{d}}{\mathrm{d}h}R (
h,\tilde{h},c ) = ( 1+h ) ^{-2} (
\frac{ ( 1+\tilde{h} ) ( c+\tilde{h} ) }{\tilde{h}}
- \frac{%
( 1+h ) ( c+h ) }{h} ) <0$
for all $h$ and $\tilde{h}$ such that $\sqrt{c}\leq\tilde{h}<h\leq H$.
Therefore, for any $H>\sqrt{c}$, there exist $\tilde{h}$ and $\delta$ such
that $\sqrt{c}<\tilde{h}\leq H$, $\delta>0$ and $R ( H,\tilde{h},c )
<-3\delta$; and thus, for sufficiently large $n$ and $p$, $%
R ( H,\tilde{h},c_{p} ) <-3\delta$. Now, $\frac{\mathrm{d}}{\mathrm{d}h}R (
h,%
\tilde{h},c_{p} ) = ( 1+h ) ^{-2} ( z_{0} ( \tilde
{h}%
) -z_{0} ( h ) ) <0$ for all $h>\tilde{h}$, as
long as $%
\tilde{h}\geq\sqrt{c_{p}}$. Hence, for sufficiently large $n$ and $p$,
$%
R ( h,\tilde{h},c_{p} ) <-3\delta$ for all $h>\tilde{h}$.
Using (%
\ref{LRasymptotica}), we get $\llvert L ( h;\mu)
\rrvert\leq e^{-{3n}\delta/{2}}O_{p} ( n^{1/2} )
=O_{p} ( e^{-n\delta} ) $ uniformly over $h\in[ H,\infty
) $.

\subsection{\texorpdfstring{Proof of Proposition \protect\ref{pr9}}{Proof of Proposition 9}}\label{sec7.5}

For brevity, we derive only the asymptotic power envelope for the
case of $\mu$-based tests. According to the
Neyman--Pearson lemma, the most powerful test of the null $\theta=0$ against
a particular alternative $\theta=\theta_{1}$ is the test which
rejects the
null when $\ln L ( \theta_{1};\mu) $ is larger than some critical
value $C$. It follows from Theorem~\ref{th7} that, for such a test to have
asymptotic size $\alpha$, $C$ must be
$C=\sqrt{V ( \theta_{1} ) }\Phi^{-1} ( 1-\alpha)
+m ( \theta_{1} )$,
where $m ( \theta_{1} ) = ( -\theta_{1}^{2}+1-e^{-\theta
_{1}^{2}} ) /4$ and $V ( \theta_{1} ) = ( \theta
_{1}^{2}-1+e^{-\theta_{1}^{2}} ) /2$ are obtained from (\ref{meanmu})
and (\ref{covariancemu}) by the re-parametrization $\theta=\sqrt{-\ln
( 1-h^{2}/c ) }$. Now, according to Le Cam's third lemma and
Theorem~\ref{th7}, under $\theta=\theta_{1}$, $\ln L ( \theta
_{1};\mu)
\stackrel{d}{\rightarrow} N ( m ( \theta_{1} ) + V
( \theta
_{1} ), V ( \theta_{1} ) ) $. Therefore, the asymptotic
power $\beta( \theta_{1};\mu) $ of the asymptotically most
powerful test of $\theta= 0$ against $\theta= \theta_{1}$ is
(\ref{localpowermu}).

\subsection{\texorpdfstring{Proof of Proposition \protect\ref{pr10}}{Proof of Proposition 10}}\label{sec7.6}

As shown by \citet{BaiBenPec05} in the complex case and by
\citet{FerPec09} in the real case, the convergence (\ref{TW}) takes place
not only under the null, but also under alternatives $h=h_{1}$ with
$h_{1}<%
\sqrt{c}$, yielding $\theta=\theta_{1}<\infty$ under the
parametrization $%
\theta=\sqrt{-\ln( 1-h^{2}/c ) }$. Hence, (\ref{powerTW}) follows.

Formulas (\ref{betaJ}) and (\ref{betaCCLR}) can be established using
conceptually similar steps. To save space, below we only establish
formula (\ref{betaCCLR}). The following technical lemma is proven in the
Supplementary Appendix [\citet{OnaMorHal}].

\begin{lemma}\label{le13}
Let $\mathrm{CLR}$ be the ``corrected'' likelihood
ratio statistic as defined in Example~\ref{ex3}. Then, under the null, as $%
p,n\rightarrow_{c} \infty$, the vector $(\mathrm{CLR},\break\Delta_{p} ( z_{0}(h) ) )
$ converges in distribution to a Gaussian vector $ (
\zeta_{1},\zeta_{2} ) $ with\break $\operatorname{Cov}(\zeta_{1},\zeta_{2} )
=-2h+2\ln( 1+h ) $.
\end{lemma}

Lemma~\ref{le13} and (\ref{equivalence2}) imply the convergence in
distribution of the vector $ ( \mathrm{CLR},\ln L ( h;\lambda) ) $ to a
Gaussian vector $ ( \zeta_{1},-\frac{1}{2}\zeta_{2} )$. From
\citet{Baietal09}, we know that, under the null, $\mathrm{CLR}
\stackrel{d}{\rightarrow
} %
N ( -\frac{1}{2}\ln( 1 - c ),-2\ln( 1 - c ) - 2c ) $. By Le Cam's third
lemma, under the alternative $h=h_{1}$, $\mathrm{CLR}$ converges to a Gaussian
random variable with the same variance but with mean equal to
$-\frac{1}{2}\ln( 1-c ) +
\operatorname{Cov}(\zeta_{1},-\frac{1}{2}\zeta_{2} ) =-\frac{1}{2}\ln(
1-c ) +h-\ln ( 1+h ) $ evaluated at $h=h_{1}$. Therefore, the power of
the ``corrected'' likelihood ratio test of asymptotic size $\alpha$
equals $1-\Phi( \Phi^{-1} ( 1-\alpha ) -\frac{h_{1}-\ln( 1+h_{1} )
}{\sqrt{-2\ln( 1-c ) -2c}} ) $. Using the reparametrization $\theta
_{1}=\sqrt{%
-\ln( 1-h_{1}^{2}/c ) }$, we get (\ref{betaCCLR}).
\end{appendix}

\section*{Acknowledgments}

This work started when the first two authors worked at and the third
author visited Columbia University. We would like to thank Tony Cai,
the Associate Editor, Nick Patterson and an anonymous referee for
helpful and encouraging comments.

\begin{supplement}
\stitle{Supplementary Appendix}
\slink[doi]{10.1214/13-AOS1100SUPP} 
\sdatatype{.pdf}
\sfilename{aos1100\_supp.pdf}
\sdescription{The Supplementary Appendix contains proofs of Lemmas~\ref{le4},
\ref{le5},~\ref{le6},~\ref{le11},~\ref{le12} and~\ref{le13}.}
\end{supplement}


\printaddresses


\begin{thebibliography}{53}

\bibitem[\protect\citeauthoryear{Bai}{1993}]{Bai93}
\begin{barticle}[mr]
\bauthor{\bsnm{Bai},~\bfnm{Z.~D.}\binits{Z.~D.}}
(\byear{1993}).
\btitle{Convergence rate of expected spectral distributions of large random
  matrices. {II}. {S}ample covariance matrices}.
\bjournal{Ann. Probab.}
\bvolume{21}
\bpages{649--672}.
\bid{issn={0091-1798}, mr={1217560}}
\bptok{imsref}%
\end{barticle}
\endbibitem

\bibitem[\protect\citeauthoryear{Bai and Silverstein}{2004}]{BaiSil04}
\begin{barticle}[mr]
\bauthor{\bsnm{Bai},~\bfnm{Z.~D.}\binits{Z.~D.}} \AND
  \bauthor{\bsnm{Silverstein},~\bfnm{Jack~W.}\binits{J.~W.}}
(\byear{2004}).
\btitle{C{LT} for linear spectral statistics of large-dimensional sample
  covariance matrices}.
\bjournal{Ann. Probab.}
\bvolume{32}
\bpages{553--605}.
\bid{doi={10.1214/aop/1078415845}, issn={0091-1798}, mr={2040792}}
\bptok{imsref}%
\end{barticle}
\endbibitem

\bibitem[\protect\citeauthoryear{Bai et~al.}{2009}]{Baietal09}
\begin{barticle}[mr]
\bauthor{\bsnm{Bai},~\bfnm{Zhidong}\binits{Z.}},
  \bauthor{\bsnm{Jiang},~\bfnm{Dandan}\binits{D.}},
  \bauthor{\bsnm{Yao},~\bfnm{Jian-Feng}\binits{J.-F.}} \AND
  \bauthor{\bsnm{Zheng},~\bfnm{Shurong}\binits{S.}}
(\byear{2009}).
\btitle{Corrections to {LRT} on large-dimensional covariance matrix by {RMT}}.
\bjournal{Ann. Statist.}
\bvolume{37}
\bpages{3822--3840}.
\bid{doi={10.1214/09-AOS694}, issn={0090-5364}, mr={2572444}}
\bptok{imsref}%
\end{barticle}
\endbibitem

\bibitem[\protect\citeauthoryear{Baik, Ben~Arous and P{\'e}ch{\'e}}{2005}]{BaiBenPec05}
\begin{barticle}[mr]
\bauthor{\bsnm{Baik},~\bfnm{Jinho}\binits{J.}},
  \bauthor{\bsnm{Ben~Arous},~\bfnm{G{\'e}rard}\binits{G.}} \AND
  \bauthor{\bsnm{P{\'e}ch{\'e}},~\bfnm{Sandrine}\binits{S.}}
(\byear{2005}).
\btitle{Phase transition of the largest eigenvalue for nonnull complex sample
  covariance matrices}.
\bjournal{Ann. Probab.}
\bvolume{33}
\bpages{1643--1697}.
\bid{doi={10.1214/009117905000000233}, issn={0091-1798}, mr={2165575}}
\bptok{imsref}%
\end{barticle}
\endbibitem

\bibitem[\protect\citeauthoryear{Baik and Silverstein}{2006}]{BaiSil06}
\begin{barticle}[mr]
\bauthor{\bsnm{Baik},~\bfnm{Jinho}\binits{J.}} \AND
  \bauthor{\bsnm{Silverstein},~\bfnm{Jack~W.}\binits{J.~W.}}
(\byear{2006}).
\btitle{Eigenvalues of large sample covariance matrices of spiked population
  models}.
\bjournal{J. Multivariate Anal.}
\bvolume{97}
\bpages{1382--1408}.
\bid{doi={10.1016/j.jmva.2005.08.003}, issn={0047-259X}, mr={2279680}}
\bptok{imsref}%
\end{barticle}
\endbibitem

\bibitem[\protect\citeauthoryear{Bejan}{2005}]{Bej05}
\begin{bmisc}[auto:STB|2013/04/09|15:24:00]
\bauthor{\bsnm{Bejan},~\bfnm{A.~Iu.}\binits{A.~I.}}
(\byear{2005}).
\bhowpublished{Largest eigenvalues and sample covariance matrices. Tracy--Widom
  and Painleve II; computational aspects and realization in S-plus with
  applications. Unpublished manuscript, Univ. Warwick}.
\bptok{imsref}%
\end{bmisc}
\endbibitem

\bibitem[\protect\citeauthoryear{Berthet and Rigollet}{2012}]{BerRig}
\begin{bmisc}[auto:STB|2013/04/09|15:24:00]
\bauthor{\bsnm{Berthet},~\bfnm{Q.}\binits{Q.}} \AND
  \bauthor{\bsnm{Rigollet},~\bfnm{P.}\binits{P.}}
(\byear{2012}).
\bhowpublished{Optimal detection of sparse principal components in high
  dimension. Available at arXiv:\arxivurl{1202.5070}}.
\bptok{imsref}%
\end{bmisc}
\endbibitem

\bibitem[\protect\citeauthoryear{Bianchi et~al.}{2011}]{Biaetal}
\begin{barticle}[auto:STB|2013/04/09|15:24:00]
\bauthor{\bsnm{Bianchi},~\bfnm{P.}\binits{P.}},
  \bauthor{\bsnm{Debbah},~\bfnm{M.}\binits{M.}},
  \bauthor{\bsnm{Maida},~\bfnm{M.}\binits{M.}} \AND
  \bauthor{\bsnm{Najim},~\bfnm{J.}\binits{J.}}
(\byear{2011}).
\btitle{Performance of statistical tests for single source detection
  using random matrix theory}.
\bjournal{IEEE Trans. Inform. Theory}
\bvolume{57}
\bpages{2400--2419}.
\bptok{imsref}%
\end{barticle}
\endbibitem

\bibitem[\protect\citeauthoryear{Birke and Dette}{2005}]{BirDet05}
\begin{barticle}[mr]
\bauthor{\bsnm{Birke},~\bfnm{Melanie}\binits{M.}} \AND
  \bauthor{\bsnm{Dette},~\bfnm{Holger}\binits{H.}}
(\byear{2005}).
\btitle{A note on testing the covariance matrix for large dimension}.
\bjournal{Statist. Probab. Lett.}
\bvolume{74}
\bpages{281--289}.
\bid{doi={10.1016/j.spl.2005.04.051}, issn={0167-7152}, mr={2189467}}
\bptok{imsref}%
\end{barticle}
\endbibitem

\bibitem[\protect\citeauthoryear{Bloemendal and Vir{\'a}g}{2012}]{BloVir}
\begin{bmisc}[auto:STB|2013/04/09|15:24:00]
\bauthor{\bsnm{Bloemendal},~\bfnm{A.}\binits{A.}} \AND
  \bauthor{\bsnm{Vir{\'a}g},~\bfnm{B.}\binits{B.}}
(\byear{2012}).
\bhowpublished{Limits of spiked random matrices I.
\textit{Probab. Theory Related Fields}.
To appear. DOI:\doiurl{10.1007/s00440-012-0443-2}.
Available at
  arXiv:\arxivurl{1011.1877v2}}.
\bptok{imsref}%
\end{bmisc}
\endbibitem

\bibitem[\protect\citeauthoryear{Butler and Wood}{2002}]{ButWoo02}
\begin{barticle}[mr]
\bauthor{\bsnm{Butler},~\bfnm{Ronald~W.}\binits{R.~W.}} \AND
  \bauthor{\bsnm{Wood},~\bfnm{Andrew T.~A.}\binits{A.~T.~A.}}
(\byear{2002}).
\btitle{Laplace approximations for hypergeometric functions with matrix
  argument}.
\bjournal{Ann. Statist.}
\bvolume{30}
\bpages{1155--1177}.
\bid{doi={10.1214/aos/1031689021}, issn={0090-5364}, mr={1926172}}
\bptok{imsref}%
\end{barticle}
\endbibitem

\bibitem[\protect\citeauthoryear{Chen, Zhang and Zhong}{2010}]{CheZhaZho10}
\begin{barticle}[mr]
\bauthor{\bsnm{Chen},~\bfnm{Song~Xi}\binits{S.~X.}},
  \bauthor{\bsnm{Zhang},~\bfnm{Li-Xin}\binits{L.-X.}} \AND
  \bauthor{\bsnm{Zhong},~\bfnm{Ping-Shou}\binits{P.-S.}}
(\byear{2010}).
\btitle{Tests for high-dimensional covariance matrices}.
\bjournal{J. Amer. Statist. Assoc.}
\bvolume{105}
\bpages{810--819}.
\bid{doi={10.1198/jasa.2010.tm09560}, issn={0162-1459}, mr={2724863}}
\bptok{imsref}%
\end{barticle}
\endbibitem

\bibitem[\protect\citeauthoryear{Dickey}{1983}]{Dic83}
\begin{barticle}[mr]
\bauthor{\bsnm{Dickey},~\bfnm{James~M.}\binits{J.~M.}}
(\byear{1983}).
\btitle{Multiple hypergeometric functions: Probabilistic interpretations and
  statistical uses}.
\bjournal{J. Amer. Statist. Assoc.}
\bvolume{78}
\bpages{628--637}.
\bid{issn={0162-1459}, mr={0721212}}
\bptok{imsref}%
\end{barticle}
\endbibitem

\bibitem[\protect\citeauthoryear{El~Karoui}{2007}]{ElK07}
\begin{barticle}[mr]
\bauthor{\bsnm{El~Karoui},~\bfnm{Noureddine}\binits{N.}}
(\byear{2007}).
\btitle{Tracy--{W}idom limit for the largest eigenvalue of a large class of
  complex sample covariance matrices}.
\bjournal{Ann. Probab.}
\bvolume{35}
\bpages{663--714}.
\bid{doi={10.1214/009117906000000917}, issn={0091-1798}, mr={2308592}}
\bptok{imsref}%
\end{barticle}
\endbibitem

\bibitem[\protect\citeauthoryear{Erdelyi}{1937}]{Erd37}
\begin{barticle}[auto:STB|2013/04/09|15:24:00]
\bauthor{\bsnm{Erdelyi},~\bfnm{A.}\binits{A.}}
(\byear{1937}).
\btitle{Beitrag zur theorie der konfluenten hypergeometrischen funktionen von
  mehreren veranderlichen}.
\bjournal{Sitzungsberichte, Akademie der Wissenschaften in Wien, Abteilung IIa,
  Mathematisch-Naturwissenschaftliche Klasse}
\bvolume{146}
\bpages{431--467}.
\bptok{imsref}%
\end{barticle}
\endbibitem

\bibitem[\protect\citeauthoryear{F{\'e}ral and P{\'e}ch{\'e}}{2009}]{FerPec09}
\begin{barticle}[mr]
\bauthor{\bsnm{F{\'e}ral},~\bfnm{Delphine}\binits{D.}} \AND
  \bauthor{\bsnm{P{\'e}ch{\'e}},~\bfnm{Sandrine}\binits{S.}}
(\byear{2009}).
\btitle{The largest eigenvalues of sample covariance matrices for a spiked
  population: Diagonal case}.
\bjournal{J. Math. Phys.}
\bvolume{50}
\bpages{073302, 33}.
\bid{doi={10.1063/1.3155785}, issn={0022-2488}, mr={2548630}}
\bptok{imsref}%
\end{barticle}
\endbibitem

\bibitem[\protect\citeauthoryear{Fisher, Sun and Gallagher}{2010}]{FisSunGal10}
\begin{barticle}[mr]
\bauthor{\bsnm{Fisher},~\bfnm{Thomas~J.}\binits{T.~J.}},
  \bauthor{\bsnm{Sun},~\bfnm{Xiaoqian}\binits{X.}} \AND
  \bauthor{\bsnm{Gallagher},~\bfnm{Colin~M.}\binits{C.~M.}}
(\byear{2010}).
\btitle{A new test for sphericity of the covariance matrix for high dimensional
  data}.
\bjournal{J. Multivariate Anal.}
\bvolume{101}
\bpages{2554--2570}.
\bid{doi={10.1016/j.jmva.2010.07.004}, issn={0047-259X}, mr={2719881}}
\bptok{imsref}%
\end{barticle}
\endbibitem


\bibitem[\protect\citeauthoryear{Guionnet and Ma{\"{\i}}da}{2005}]{GuiMad05}
\begin{barticle}[mr]
\bauthor{\bsnm{Guionnet},~\bfnm{A.}\binits{A.}} \AND
  \bauthor{\bsnm{Ma{\"{\i}}da},~\bfnm{M.}\binits{M.}}
(\byear{2005}).
\btitle{A {F}ourier view on the {$R$}-transform and related asymptotics of
  spherical integrals}.
\bjournal{J. Funct. Anal.}
\bvolume{222}
\bpages{435--490}.
\bid{doi={10.1016/j.jfa.2004.09.015}, issn={0022-1236}, mr={2132396}}
\bptok{imsref}%
\end{barticle}
\endbibitem

\bibitem[\protect\citeauthoryear{Hallin and Paindaveine}{2006}]{HalPai06}
\begin{barticle}[mr]
\bauthor{\bsnm{Hallin},~\bfnm{Marc}\binits{M.}} \AND
  \bauthor{\bsnm{Paindaveine},~\bfnm{Davy}\binits{D.}}
(\byear{2006}).
\btitle{Semiparametrically efficient rank-based inference for shape. {I}.
  {O}ptimal rank-based tests for sphericity}.
\bjournal{Ann. Statist.}
\bvolume{34}
\bpages{2707--2756}.
\bid{doi={10.1214/009053606000000731}, issn={0090-5364}, mr={2329465}}
\bptok{imsref}%
\end{barticle}
\endbibitem

\bibitem[\protect\citeauthoryear{Hillier}{2001}]{Hil01}
\begin{barticle}[mr]
\bauthor{\bsnm{Hillier},~\bfnm{Grant}\binits{G.}}
(\byear{2001}).
\btitle{The density of a quadratic form in a vector uniformly distributed on
  the {$n$}-sphere}.
\bjournal{Econometric Theory}
\bvolume{17}
\bpages{1--28}.
\bid{doi={10.1017/S026646660117101X}, issn={0266-4666}, mr={1863565}}
\bptok{imsref}%
\end{barticle}
\endbibitem

\bibitem[\protect\citeauthoryear{Hoyle}{2008}]{Hoy08}
\begin{barticle}[auto:STB|2013/04/09|15:24:00]
\bauthor{\bsnm{Hoyle},~\bfnm{D.~C.}\binits{D.~C.}}
(\byear{2008}).
\btitle{Automatic PCA dimension selection for high dimensional data and small
  sample sizes}.
\bjournal{J. Mach. Learn. Res.}
\bvolume{9}
\bpages{2733--2759}.
\bptok{imsref}%
\end{barticle}
\endbibitem

\bibitem[\protect\citeauthoryear{James}{1964}]{Jam64}
\begin{barticle}[mr]
\bauthor{\bsnm{James},~\bfnm{Alan~T.}\binits{A.~T.}}
(\byear{1964}).
\btitle{Distributions of matrix variates and latent roots derived from normal
  samples}.
\bjournal{Ann. Math. Statist.}
\bvolume{35}
\bpages{475--501}.
\bid{issn={0003-4851}, mr={0181057}}
\bptok{imsref}%
\end{barticle}
\endbibitem

\bibitem[\protect\citeauthoryear{John}{1971}]{Joh71}
\begin{barticle}[mr]
\bauthor{\bsnm{John},~\bfnm{S.}\binits{S.}}
(\byear{1971}).
\btitle{Some optimal multivariate tests}.
\bjournal{Biometrika}
\bvolume{58}
\bpages{123--127}.
\bid{issn={0006-3444}, mr={0275568}}
\bptok{imsref}%
\end{barticle}
\endbibitem

\bibitem[\protect\citeauthoryear{John}{1972}]{Joh72}
\begin{barticle}[mr]
\bauthor{\bsnm{John},~\bfnm{S.}\binits{S.}}
(\byear{1972}).
\btitle{The distribution of a statistic used for testing sphericity of normal
  distributions}.
\bjournal{Biometrika}
\bvolume{59}
\bpages{169--173}.
\bid{issn={0006-3444}, mr={0312619}}
\bptok{imsref}%
\end{barticle}
\endbibitem

\bibitem[\protect\citeauthoryear{Johnstone}{2001}]{Joh01}
\begin{barticle}[mr]
\bauthor{\bsnm{Johnstone},~\bfnm{Iain~M.}\binits{I.~M.}}
(\byear{2001}).
\btitle{On the distribution of the largest eigenvalue in principal components
  analysis}.
\bjournal{Ann. Statist.}
\bvolume{29}
\bpages{295--327}.
\bid{doi={10.1214/aos/1009210544}, issn={0090-5364}, mr={1863961}}
\bptok{imsref}%
\end{barticle}
\endbibitem

\bibitem[\protect\citeauthoryear{Kritchman and Nadler}{2008}]{KriNad08}
\begin{barticle}[auto:STB|2013/04/09|15:24:00]
\bauthor{\bsnm{Kritchman},~\bfnm{S.}\binits{S.}} \AND
  \bauthor{\bsnm{Nadler},~\bfnm{B.}\binits{B.}}
(\byear{2008}).
\btitle{Determining the number of components in a factor model from limited
  noisy data}.
\bjournal{Chemometrics and Intelligent Laboratory Systems}
\bvolume{94}
\bpages{19--32}.
\bptok{imsref}%
\end{barticle}
\endbibitem

\bibitem[\protect\citeauthoryear{Kritchman and Nadler}{2009}]{KriNad09}
\begin{barticle}[mr]
\bauthor{\bsnm{Kritchman},~\bfnm{Shira}\binits{S.}} \AND
  \bauthor{\bsnm{Nadler},~\bfnm{Boaz}\binits{B.}}
(\byear{2009}).
\btitle{Non-parametric detection of the number of signals: {H}ypothesis testing
  and random matrix theory}.
\bjournal{IEEE Trans. Signal Process.}
\bvolume{57}
\bpages{3930--3941}.
\bid{doi={10.1109/TSP.2009.2022897}, issn={1053-587X}, mr={2683143}}
\bptok{imsref}%
\end{barticle}
\endbibitem

\bibitem[\protect\citeauthoryear{Le~Cam}{1960}]{LeC60}
\begin{barticle}[mr]
\bauthor{\bsnm{Le~Cam},~\bfnm{Lucien}\binits{L.}}
(\byear{1960}).
\btitle{Locally asymptotically normal families of distributions. {C}ertain
  approximations to families of distributions and their use in the theory of
  estimation and testing hypotheses}.
\bjournal{Univ. California Publ. Statist.}
\bvolume{3}
\bpages{37--98}.
\bid{mr={0126903}}
\bptok{imsref}%
\end{barticle}
\endbibitem

\bibitem[\protect\citeauthoryear{Ledoit and Wolf}{2002}]{LedWol02}
\begin{barticle}[mr]
\bauthor{\bsnm{Ledoit},~\bfnm{Olivier}\binits{O.}} \AND
  \bauthor{\bsnm{Wolf},~\bfnm{Michael}\binits{M.}}
(\byear{2002}).
\btitle{Some hypothesis tests for the covariance matrix when the dimension is
  large compared to the sample size}.
\bjournal{Ann. Statist.}
\bvolume{30}
\bpages{1081--1102}.
\bid{doi={10.1214/aos/1031689018}, issn={0090-5364}, mr={1926169}}
\bptok{imsref}%
\end{barticle}
\endbibitem

\bibitem[\protect\citeauthoryear{Lijoi and Regazzini}{2004}]{LijReg04}
\begin{barticle}[mr]
\bauthor{\bsnm{Lijoi},~\bfnm{Antonio}\binits{A.}} \AND
  \bauthor{\bsnm{Regazzini},~\bfnm{Eugenio}\binits{E.}}
(\byear{2004}).
\btitle{Means of a {D}irichlet process and multiple hypergeometric functions}.
\bjournal{Ann. Probab.}
\bvolume{32}
\bpages{1469--1495}.
\bid{doi={10.1214/009117904000000270}, issn={0091-1798}, mr={2060305}}
\bptok{imsref}%
\end{barticle}
\endbibitem

\bibitem[\protect\citeauthoryear{Mauchly}{1940}]{Mau40}
\begin{barticle}[mr]
\bauthor{\bsnm{Mauchly},~\bfnm{John~W.}\binits{J.~W.}}
(\byear{1940}).
\btitle{Significance test for sphericity of a normal {$n$}-variate
  distribution}.
\bjournal{Ann. Math. Statist.}
\bvolume{11}
\bpages{204--209}.
\bid{issn={0003-4851}, mr={0002084}}
\bptok{imsref}%
\end{barticle}
\endbibitem

\bibitem[\protect\citeauthoryear{Mo}{2012}]{Mo12}
\begin{barticle}[mr]
\bauthor{\bsnm{Mo},~\bfnm{M.~Y.}\binits{M.~Y.}}
(\byear{2012}).
\btitle{Rank 1 real {W}ishart spiked model}.
\bjournal{Comm. Pure Appl. Math.}
\bvolume{65}
\bpages{1528--1638}.
\bid{doi={10.1002/cpa.21415}, issn={0010-3640}, mr={2969495}}
\bptnote{check year}%
\bptok{imsref}%
\end{barticle}
\endbibitem

\bibitem[\protect\citeauthoryear{Nadakuditi and Edelman}{2008}]{NadEde08}
\begin{barticle}[mr]
\bauthor{\bsnm{Nadakuditi},~\bfnm{Raj~Rao}\binits{R.~R.}} \AND
  \bauthor{\bsnm{Edelman},~\bfnm{Alan}\binits{A.}}
(\byear{2008}).
\btitle{Sample eigenvalue based detection of high-dimensional signals in white
  noise using relatively few samples}.
\bjournal{IEEE Trans. Signal Process.}
\bvolume{56}
\bpages{2625--2638}.
\bid{doi={10.1109/TSP.2008.917356}, issn={1053-587X}, mr={1500236}}
\bptok{imsref}%
\end{barticle}
\endbibitem

\bibitem[\protect\citeauthoryear{Nadakuditi and Silverstein}{2010}]{NadSil10}
\begin{barticle}[auto:STB|2013/04/09|15:24:00]
\bauthor{\bsnm{Nadakuditi},~\bfnm{R.~R.}\binits{R.~R.}} \AND
  \bauthor{\bsnm{Silverstein},~\bfnm{J.~W.}\binits{J.~W.}}
(\byear{2010}).
\btitle{Fundamental limit of sample generalized eigenvalue based detection of
  signals in noise using relatively few signal-bearing and noise-only samples}.
\bjournal{IEEE Journal of Selected Topics in Signal Processing}
\bvolume{4}
\bpages{468--480}.
\bptok{imsref}%
\end{barticle}
\endbibitem

\bibitem[\protect\citeauthoryear{Nadler}{2008}]{Nad08}
\begin{barticle}[mr]
\bauthor{\bsnm{Nadler},~\bfnm{Boaz}\binits{B.}}
(\byear{2008}).
\btitle{Finite sample approximation results for principal component analysis: A
  matrix perturbation approach}.
\bjournal{Ann. Statist.}
\bvolume{36}
\bpages{2791--2817}.
\bid{doi={10.1214/08-AOS618}, issn={0090-5364}, mr={2485013}}
\bptok{imsref}%
\end{barticle}
\endbibitem

\bibitem[\protect\citeauthoryear{Olver}{1997}]{Olv97}
\begin{bbook}[mr]
\bauthor{\bsnm{Olver},~\bfnm{Frank W.~J.}\binits{F.~W.~J.}}
(\byear{1997}).
\btitle{Asymptotics and Special Functions}.
\bpublisher{AK Peters}, \blocation{Wellesley, MA}.
\bid{mr={1429619}}
\bptok{imsref}%
\end{bbook}
\endbibitem

\bibitem[\protect\citeauthoryear{Onatski}{2009}]{Ona09}
\begin{barticle}[mr]
\bauthor{\bsnm{Onatski},~\bfnm{Alexei}\binits{A.}}
(\byear{2009}).
\btitle{Testing hypotheses about the numbers of factors in large factor
  models}.
\bjournal{Econometrica}
\bvolume{77}
\bpages{1447--1479}.
\bid{doi={10.3982/ECTA6964}, issn={0012-9682}, mr={2561070}}
\bptok{imsref}%
\end{barticle}
\endbibitem

\bibitem[\protect\citeauthoryear{Onatski}{2010}]{Ona10}
\begin{barticle}[auto:STB|2013/04/09|15:24:00]
\bauthor{\bsnm{Onatski},~\bfnm{A.}\binits{A.}}
(\byear{2010}).
\btitle{Determining the number of factors from empirical distribution of
  eigenvalues}.
\bjournal{Rev. Econom. Statist.}
\bvolume{92}
\bpages{1004--1016}.
\bptok{imsref}%
\end{barticle}
\endbibitem

\bibitem[\protect\citeauthoryear{Onatski, Moreira and Hallin}{2013}]{OnaMorHal}
\begin{bmisc}[auto:STB|2013/04/09|15:24:00]
\bauthor{\bsnm{Onatski},~\bfnm{A.}\binits{A.}},
  \bauthor{\bsnm{Moreira},~\bfnm{M.~J.}\binits{M.~J.}} \AND
  \bauthor{\bsnm{Hallin},~\bfnm{M.}\binits{M.}}
(\byear{2013}).
\bhowpublished{Supplement to ``Asymptotic power of sphericity tests for
  high-dimensional data.'' DOI:\doiurl{10.1214/13-AOS1100SUPP}}.
\bptok{imsref}%
\end{bmisc}
\endbibitem

\bibitem[\protect\citeauthoryear{Patterson, Price and Reich}{2006}]{PatPriRei06}
\begin{barticle}[auto:STB|2013/04/09|15:24:00]
\bauthor{\bsnm{Patterson},~\bfnm{N.}\binits{N.}},
  \bauthor{\bsnm{Price},~\bfnm{A.~L.}\binits{A.~L.}} \AND
  \bauthor{\bsnm{Reich},~\bfnm{D.}\binits{D.}}
(\byear{2006}).
\btitle{Population structure and eigenanalysis}.
\bjournal{PLoS Genetics}
\bvolume{2}
\bpages{2074--2093}.
\bptok{imsref}%
\end{barticle}
\endbibitem

\bibitem[\protect\citeauthoryear{Perry and Wolfe}{2010}]{PerWol10}
\begin{barticle}[auto:STB|2013/04/09|15:24:00]
\bauthor{\bsnm{Perry},~\bfnm{P.~O.}\binits{P.~O.}} \AND
  \bauthor{\bsnm{Wolfe},~\bfnm{P.~J.}\binits{P.~J.}}
(\byear{2010}).
\btitle{Minimax rank estimation for subspace tracking}.
\bjournal{IEEE Journal of Selected Topics in Signal Processing}
\bvolume{4}
\bpages{504--513}.
\bptok{imsref}%
\end{barticle}
\endbibitem



\bibitem[\protect\citeauthoryear{Schott}{2006}]{Sch06}
\begin{barticle}[mr]
\bauthor{\bsnm{Schott},~\bfnm{James~R.}\binits{J.~R.}}
(\byear{2006}).
\btitle{A high-dimensional test for the equality of the smallest eigenvalues of
  a covariance matrix}.
\bjournal{J. Multivariate Anal.}
\bvolume{97}
\bpages{827--843}.
\bid{doi={10.1016/j.jmva.2005.05.003}, issn={0047-259X}, mr={2256563}}
\bptok{imsref}%
\end{barticle}
\endbibitem

\bibitem[\protect\citeauthoryear{Silverstein and Bai}{1995}]{SilBai95}
\begin{barticle}[mr]
\bauthor{\bsnm{Silverstein},~\bfnm{Jack~W.}\binits{J.~W.}} \AND
  \bauthor{\bsnm{Bai},~\bfnm{Z.~D.}\binits{Z.~D.}}
(\byear{1995}).
\btitle{On the empirical distribution of eigenvalues of a class of
  large-dimensional random matrices}.
\bjournal{J. Multivariate Anal.}
\bvolume{54}
\bpages{175--192}.
\bid{doi={10.1006/jmva.1995.1051}, issn={0047-259X}, mr={1345534}}
\bptok{imsref}%
\end{barticle}
\endbibitem

\bibitem[\protect\citeauthoryear{Srivastava}{2005}]{Sri05}
\begin{barticle}[mr]
\bauthor{\bsnm{Srivastava},~\bfnm{Muni~S.}\binits{M.~S.}}
(\byear{2005}).
\btitle{Some tests concerning the covariance matrix in high dimensional data}.
\bjournal{J. Japan Statist. Soc.}
\bvolume{35}
\bpages{251--272}.
\bid{issn={1882-2754}, mr={2328427}}
\bptok{imsref}%
\end{barticle}
\endbibitem

\bibitem[\protect\citeauthoryear{Srivastava and Karlsson}{1985}]{SriKar85}
\begin{bbook}[mr]
\bauthor{\bsnm{Srivastava},~\bfnm{H.~M.}\binits{H.~M.}} \AND
  \bauthor{\bsnm{Karlsson},~\bfnm{Per~W.}\binits{P.~W.}}
(\byear{1985}).
\btitle{Multiple {G}aussian Hypergeometric Series}.
\bpublisher{Ellis Horwood}, \blocation{Chichester}.
\bid{mr={0834385}}
\bptok{imsref}%
\end{bbook}
\endbibitem

\bibitem[\protect\citeauthoryear{Sugiura}{1972}]{Sug72}
\begin{barticle}[mr]
\bauthor{\bsnm{Sugiura},~\bfnm{Nariaki}\binits{N.}}
(\byear{1972}).
\btitle{Locally best invariant test for sphericity and the limiting
  distributions}.
\bjournal{Ann. Math. Statist.}
\bvolume{43}
\bpages{1312--1316}.
\bid{issn={0003-4851}, mr={0311032}}
\bptok{imsref}%
\end{barticle}
\endbibitem

\bibitem[\protect\citeauthoryear{Uhlig}{1994}]{Uhl94}
\begin{barticle}[mr]
\bauthor{\bsnm{Uhlig},~\bfnm{Harald}\binits{H.}}
(\byear{1994}).
\btitle{On singular {W}ishart and singular multivariate beta distributions}.
\bjournal{Ann. Statist.}
\bvolume{22}
\bpages{395--405}.
\bid{doi={10.1214/aos/1176325375}, issn={0090-5364}, mr={1272090}}
\bptok{imsref}%
\end{barticle}
\endbibitem

\bibitem[\protect\citeauthoryear{van~der Vaart}{1998}]{van98}
\begin{bbook}[mr]
\bauthor{\bparticle{van~der} \bsnm{Vaart},~\bfnm{A.~W.}\binits{A.~W.}}
(\byear{1998}).
\btitle{Asymptotic Statistics}.
\bseries{Cambridge Series in Statistical and Probabilistic Mathematics}
\bvolume{3}.
\bpublisher{Cambridge Univ. Press}, \blocation{Cambridge}.
\bid{mr={1652247}}
\bptok{imsref}%
\end{bbook}
\endbibitem

\bibitem[\protect\citeauthoryear{Wang}{2012}]{Wan12}
\begin{barticle}[mr]
\bauthor{\bsnm{Wang},~\bfnm{Dong}\binits{D.}}
(\byear{2012}).
\btitle{The largest eigenvalue of real symmetric, {H}ermitian and {H}ermitian
  self-dual random matrix models with rank one external source, {P}art {I}}.
\bjournal{J. Stat. Phys.}
\bvolume{146}
\bpages{719--761}.
\bid{doi={10.1007/s10955-012-0417-x}, issn={0022-4715}, mr={2916094}}
\bptok{imsref}%
\end{barticle}
\endbibitem

\end{thebibliography}
\end{document}